# ANALYTIC HYPOELLIPTICITY IN DIMENSION TWO

MICHAEL CHRIST

ABSTRACT. A simple geometric condition is sufficient for analytic hypoellipticity of sums of squares of two vector fields in $\mathbb{R}^2$. This condition is proved to be necessary for generic vector fields and for various special cases, and to be both necessary and sufficient for a closely related family of operators.

## CONTENTS

## 1. INTRODUCTION

Consider a linear partial differential operator $L$, representable as a sum of squares $\sum_j X_j^2$ of finitely many real vector fields with real analytic coefficients. Under what conditions is $L$ analytic hypoelliptic? Many examples and some partial results are known, but the techniques presently available fall far short of resolving the question, and no plausible candidate for a necessary and sufficient condition has been put forward. In this paper a nearly complete treatment will be given of the simplest possible case, that of two vector fields in $\mathbb{R}^2$. The main results are:

*Date*: February 4, 1996.

Research supported by NSF grant DMS-9306833 and at MSRI by NSF grant DMS-9022140.

- Formulation of a conjectured necessary and sufficient condition for analytic hypoellipticity.
- Proof of sufficiency of this condition.
- Complete reduction of the question of necessity to an eigenvalue problem for certain ordinary differential operators.
- Solution of generic eigenvalue problems, and hence proof of the conjecture for generic vector fields.
- Treatment of certain examples.
- Formulation and proof of a necessary and sufficient condition for microlocal analytic hypoellipticity of $(X_1 + iX_2) \circ (X_1 - iX_2)$, under a natural pseudoconvexity hypothesis.
- Introduction of a new geometric invariant, a rational number $q$, associated to a pair of vector fields in $\mathbb{R}^2$.

One upshot of this analysis is that analytic hypoellipticity is very rare.

**Definition.** Let $\{X_1, X_2\}, \{Y_1, Y_2\}$ be two pairs of real vector fields with real analytic coefficients, defined in a neighborhood of a point $p$. We say that

$$\mathrm{span}\{X_1, X_2\} \equiv \mathrm{span}\{Y_1, Y_2\}$$

in a neighborhood $U$ of $p$ if each $X_j$ may be expressed as a linear combination, with analytic coefficients, of the $Y_i$, and conversely each $Y_i$ may be so expressed in terms of the $X_j$, in $U$.

Let $X_1, X_2$ be real vector fields with $C^\omega$ coefficients, defined in a neighborhood of $p \in \mathbb{R}^2$. Suppose that the Lie algebra generated by them spans the tangent space at $p$.

**Conjecture 1.1.** *$L$ is analytic hypoelliptic in some neighborhood of $p$ if and only if there exists a system of coordinates $(x, t)$ with origin at $p$, in which*

$$\mathrm{span}\{X_1, X_2\} \equiv \mathrm{span}\{\partial_x, x^{m-1}\partial_t\} \tag{1.1}$$

*in some neighborhood of the origin, for some $m \geq 1$.*

That (1.1) implies analytic hypoellipticity is very well known in the elliptic case $m = 1$ and symplectic case $m = 2$. If a point $p_0$ is fixed, then within the class of all pairs of $C^\omega$ vector fields for which $X_1, X_2, [X_1, X_2]$ are linearly dependent at $p_0$ but which satisfy the bracket hypothesis at $p_0$ to some fixed order, the condition of Conjecture 1.1 is violated generically. On the other hand, for any pair of vector fields satisfying the bracket hypothesis, the set of points $p \in \mathbb{R}^2$ at which the condition of Conjecture 1.1 fails to hold is discrete.

**Theorem 1.2.** *Suppose that for some $m \geq 1$ there exist coordinates $(x, t)$ in a neighborhood $U$ of $0$ in which $\operatorname{span}\{X_1, X_2\} \equiv \operatorname{span}\{\partial_x, x^{m-1}\partial_t\}$. Then $X_1^2 + X_2^2$ is analytic hypoelliptic in $U$.*

This is a special case of a relatively elementary theorem of Grušin [12]. We include a proof in Section 7 in order to show how it follows from the same point of view, and from a subset of the machinery, that we will develop in demonstrating nonhypoellipticity in other cases. With this result in hand, our aim is to prove that analytic hypoellipticity fails in all nonelementary cases.

A few definitions are required before our main results can be properly formulated. The Lie algebra generated by $X_1, X_2$ is always assumed to span the tangent space. By a normalized homogeneous polynomial of degree $m - 1$ we will mean a function $Q(x, z)$ of the form

$$Q(x, z) = x^{m-1} + \sum_{j=0}^{m-3} c_j z^{m-1-j} x^j, \tag{1.2}$$

for some $m \geq 1$, where each $c_j \in \mathbb{R}$.[1]

**Definition.** Two normalized homogeneous polynomials $Q, P$ of degree $m - 1$ are equivalent if there exists a nonzero real constant $c$ such that for all $(x, t)$,

$$P(x, t) = Q(x, ct).$$

Define $m = m(p)$ to be the type of the set $\{X_1, X_2\}$ at $p$. This invariant is 1 if the vector fields are linearly independent at $p$, and otherwise is defined to be the smallest integer such that $X_1, X_2$ together with all of their iterated Lie brackets having $m$ or fewer factors span the tangent space to $\mathbb{R}^2$ at $p$. Denote by $\mathbb{Q}^+$ the set of all positive rational numbers.

The following lemma will be proved in Section 2.

**Lemma 1.3.** *Suppose that $X_1, X_2$ are real vector fields with $C^\omega$ coefficients defined near $p \in \mathbb{R}^2$. Let $m \geq 1$ be their type at $p$. Suppose that there do not exist coordinates $(x, t)$ with origin at $p$ such that $\operatorname{span}\{X_1, X_2\} \equiv \operatorname{span}\{\partial_x, x^{m-1}\partial_t\}$ in a neighborhood of $0$. Then there exist $q \in \mathbb{Q}^+$, coordinates $(x, t)$ with origin at $p$, and an analytic function $\tilde{\Theta}(x, t)$ of the form*

$$\tilde{\Theta}(x, t) = x^{m-1} + \sum_{j=0}^{m-3} \beta_j(t) x^j$$

---

[1] If $m < 3$ then the sum over $j$ is taken to be empty.

*such that near* 0,

$$\text{span}\{X_1, X_2\} \equiv \text{span}\{\partial_x, \tilde{\Theta}(x, t)\partial_t\}$$

*with the following properties: Each* $\beta_j \in C^\omega$ *is real-valued,* $\beta_j(0) = 0$, *each* $\beta_j(t) = c_j t^{(m-j-1)q}(1 + O(|t|^\delta))$ *for some* $\delta > 0$, *and at least one coefficient* $c_j$ *is nonzero.*

*The quantity* $q$ *is independent of the coordinate system used and choices made in its definition. The polynomial* $Q(x, z) = x^{m-1} + \sum_{j=0}^{m-3} c_j z^{m-j-1} x^j$ *is likewise invariant, modulo the equivalence relation defined above.*

Note that necessarily $c_j = 0$ whenever the exponent $(m - j - 1)q$ is not an integer.

To a pair $(L, p)$, where $L = X_1^2 + X_2^2$ and $p \in \mathbb{R}^2$, we associate a family of ordinary differential operators

$$\mathcal{L}_z = -\partial_x^2 + Q^2(x, z), \qquad (x, z) \in \mathbb{R} \times \mathbb{C}$$

depending polynomially on the complex parameter $z$, where $Q$ is a member of the equivalence class of normalized homogeneous polynomials associated to $L$ via Lemma 1.3.

**Definition.**

$$NE(L, p) = \{z \in \mathbb{C} : \text{ there exists } 0 \neq f \in \mathcal{S}(\mathbb{R}) \text{ satisfying } \mathcal{L}_z f \equiv 0\}.$$
$$\text{(1.3)}$$

This definition does not quite make sense, because associated to $L$ is not a single family of ordinary differential operators $\mathcal{L}_z$, but rather an equivalence class. Therefore we define two subsets $S, S'$ of $\mathbb{C}$ to be equivalent if there exists a nonzero real number $\gamma$ such that $S' = \gamma S = \{\gamma \zeta : \zeta \in S\}$. $NE(L, p)$ is then defined to be the equivalence class of the set described in (1.3).

When we are discussing a family of ordinary differential operators $\mathcal{L}_z$ of the form (1.3), not necessarily associated to any partial differential operator, we denote the same set of nonlinear eigenvalues by $E(Q)$.

It is a general principle that, at least for certain classes of partial differential operators, the question of analytic hypoellipticity is linked to global eigenvalue problems. In the present context this link may be formulated rigorously, as follows.

**Theorem 1.4.** *Suppose that* $L$ *is a sum of squares of two real,* $C^\omega$ *vector fields in* $\mathbb{R}^2$ *satisfying the bracket hypothesis. If* $NE(L, p)$ *is nonempty, then* $L$ *is not analytic hypoelliptic in any neighborhood of* $p$.

**Conjecture 1.5.** *If* $Q$ *is a normalized homogeneous polynomial of degree* $m - 1 \geq 1$, *then* $E(Q) = \emptyset$ *if and only if* $Q(x, z) \equiv x^{m-1}$.

Three pieces of evidence will be offered to support this conjecture. Firstly, in Section 5 the analogous assertion for a family of ordinary differential operators that is formally very similar to $\mathcal{L}_z$ will be formulated and proved. Secondly, representative special cases, as listed in Proposition 1.7, will be treated. Thirdly, the conjecture will be proved to be valid for generic $Q$.

The following notion of genericity is natural in this discussion. Recall that a subset $E$ of $\mathbb{C}^d$ is said to be pluripolar if there exists a nonconstant plurisubharmonic function $h : \mathbb{C}^d \mapsto [-\infty, \infty]$, satisfying $h \equiv -\infty$ on $E$. A subset of $\mathbb{R}^d$ is said to be pluripolar if it is pluripolar as a subset of $\mathbb{C}^d$. We identify the set of all normalized homogeneous polynomials $Q$ of degree $m-1$ with $\mathbb{R}^{m-2} \subset \mathbb{C}^{m-2}$ via the correspondence $\zeta \mapsto Q_\zeta(x, z) = x^{m-1} + \sum_{j=0}^{m-3} \zeta_{j+1} z^{m-1-j} x^j$. In Section 5 it will be shown that Conjecture 1.5 is valid for generic $Q$, in this sense. This and Theorem 1.4 together have the following consequence.

**Theorem 1.6.** *For each $m \geq 3$, there exists a pluripolar set $B_m \subset \mathbb{R}^{m-2}$ such that whenever the pair $\{X_1, X_2\}$ is of type $m$ at $p$ and some normalized homogeneous polynomial $Q$ associated to $L$ belongs to $\mathbb{R}^{m-2} \backslash B_m$, $L$ fails to be analytic hypoelliptic in any neighborhood of $p$.*

The following is a list of examples for which it is relatively easy to show that $E(Q)$ is indeed nonempty.

**Proposition 1.7.** *Let $Q(x, z) = x^{m-1} - z^k x^{m-k-1}$. If at least one of the following conditions holds, then $E(Q) \neq \emptyset$, and consequently for any positive integer $r$ and for either choice of the $\pm$ sign,*

$$L = \partial_x^2 + \left([x^{m-1} \pm t^r x^{m-k-1}]\partial_t\right)^2$$

*is not analytic hypoelliptic in any neighborhood of $0$ :*

- *$m/k$ is not an integer*
- *$m$ is odd*
- *$k$ is even*
- *$m$ is divisible by $4$.*

We turn now to the analysis of the closely related operators

$$L = (X_1 + iX_2)(X_1 - iX_2) + c_1 X_1 + c_2 X_2 + c_3$$

where the $c_j$ are $C^\omega$, complex valued coefficients; the lower order terms $c_1 X_1 + c_2 X_2 + c_3$ play no essential role. These are lower-dimensional analogues of the operators $\bar{\partial}_b \bar{\partial}_b^*$ which arise on three-dimensional CR manifolds, and it is natural to impose a supplemental hypothesis.

**Definition.** $\{X_1, X_2\}$ is a pseudoconvex pair in a neighborhood of $p \in \mathbb{R}^2$ if either $X_1, X_2$ span $\mathbb{R}^2$ at $p$, or they are dependent at $p$ and there exist a real, smooth vector field $T$ transverse to span$\{X_1, X_2\}$ at $p$, and smooth coefficients $h, b_1, b_2$ such that $[X_1, X_2] \equiv hT + b_1 X_1 + b_2 X_2$ and $h$ does not change sign in some neighborhood of $p$.

This condition is clearly independent of $T$ and of the choice of basis for span$\{X_1, X_2\}$. However, choosing $T$ does not uniquely determine the coefficients $h, b_j$.

If $X_1, X_2$ are linearly independent at $p$ then $L$ is elliptic, hence analytic hypoelliptic, near $p$. In the linearly dependent case, assuming the bracket hypothesis to be satisfied, let $(x, t)$ be a coordinate system with the properties of Lemma 1.3. Let $(\xi, \tau)$ be Fourier variables dual to $(x, t)$. The characteristic variety $\Sigma$ of $L$ is near $0 \in \mathbb{R}^2$ a trivial line bundle over the analytic variety $V = \{(x, t) \in \mathbb{R}^2 : \tilde{\Theta}(x, t) = 0\}$ consisting of all points in the base space at which $X_1, X_2$ are dependent [2]. Restricting attention to a neighborhood of 0 in which $V$ is connected, $\Sigma$ splits in a unique way as a union of two half line bundles $\Sigma^{\pm}$, where each half line has vertex at $\xi = \tau = 0$. It will be shown in Section 8 that exactly one of these two, denoted $\Sigma^+$, has a conic neighborhood $\Gamma^+$ in which the principal symbols $\sigma_1$ of the vector fields satisfy

$$\sigma_1(i[X_1, X_2]) = \mu + a_1 \sigma_1(iX_1) + a_2 \sigma_1(iX_2)$$

for some functions $\mu, a_1, a_2$ such that $\mu \leq 0$. From this it is straightforward to deduce, following the method of Kohn [14] and exploiting the bracket hypothesis and Garding's inequality, that $L$ is microlocally $C^\infty$ hypoelliptic in $\Gamma^+$.

Our final result characterizes microlocal analytic hypoellipticity of $L$ in $\Gamma^+$. Denote by $WF_a(u)$ the analytic wave front set, as defined in [17].

**Theorem 1.8.** *Assume $X_1, X_2$ are linearly dependent at $p$. Under the hypotheses of finite type and pseudoconvexity,*

$$WF_a(u) \cap \Gamma^+ \subset WF_a(Lu) \cap \Gamma^+ \text{ for all distributions } u$$

*if and only if there exist $m \geq 2$ and coordinates $(x, t)$ with origin at $p$ such that*

$$\text{span}\{X_1, X_2\} \equiv \text{span}\{\partial_x, x^{m-1}\partial_t\}$$

*in a neighborhood of 0.*

Assuming pseudoconvexity, $L$ will not be microlocally analytic hypoelliptic in any conic neighborhood of $\Sigma^-$, at least in the special case where $L = (X_1 + iX_2)(X_1 - iX_2)$. Indeed, the theorem of Trepreau [22] asserts that $(X_1 - iX_2)$ is not microlocally analytic hypoelliptic there, hence the composition cannot be. The proof of Theorem 1.8 will be outlined in Section 8.

---

[2] The hypotheses of pseudoconvexity and finite type force $V$ to have real dimension one.

The principal results of this paper were announced in [4]. The invariant $q$ is already ubiquitous in our estimates, but in terms of Gevrey class hypoellipticity has a deeper significance outlined in [4], which will be the subject of a subsequent publication.

The technique used here to reduce analytic (non)hypoellipticity to eigenvalue problems for ordinary differential equations is identical to that used by the author in establishing a counterexample to global analytic regularity [3]; only the details of the formulas must be changed. Furthermore, various lemmas are contained, at least implicitly, in earlier works [5],[6],[8]. Therefore this paper is not entirely self-contained.

## 2. Coordinates and canonical forms

In this section we prove Lemma 1.3. Assume that $X_1, X_2$ are linearly dependent at $p \in \mathbb{R}^2$. The Lie bracket hypothesis ensures that at least one of $X_1, X_2$ does not vanish at $p$, so there exists a $C^\omega$ coordinate system $(x, t)$ with origin at $p$ in which $\text{span}\{X_1, X_2\} \equiv \text{span}\{\partial_x, a(x,t)\partial_t\}$ for some $a \in C^\omega$. By the Weierstraß preparation theorem, this span is in turn equal to the span of $\partial_x$ and $\tilde{\Theta}(x,t)\partial_t$, for some $\tilde{\Theta} \in C^\omega$ of the form

$$\tilde{\Theta}(x,t) = x^{m-1} + \sum_{j=0}^{m-2} \beta_j(t) x^j$$

for some $m \geq 1$, where each $\beta_j \in C^\omega$ is real-valued and $\beta_j(0) = 0$. Matters may be further reduced to the case where $\beta_{m-2} \equiv 0$ by a change of variables $(x,t) \mapsto (x - (m-1)^{-1}\beta_{m-2}(t), t)$, and we shall always do so. Thus

$$\tilde{\Theta}(x,t) = x^{m-1} + \sum_{j=0}^{m-3} \beta_j(t) x^j. \tag{2.1}$$

From the expression for $\tilde{\Theta}$ one finds that the quantity $m$ defined in this way equals the type as defined in Section 1.

Define $\tau_j$ to be the order of vanishing of $\beta_j$ at $t = 0$, that is, $\beta_j(t) = d_j t^{\tau_j} + O(t^{\tau_j+1})$ for some $d_j \neq 0$.

**Definition.**

$$q = \begin{cases} \infty & \text{if each } \beta_j \text{ vanishes identically} \\ \min_j \tau_j/(m-1-j) & \text{otherwise.} \end{cases}$$

The lowest order part of $\tilde{\Theta}$ will be denoted by $\Theta$.

**Definition.** If $q < \infty$,

$$\Theta(x,t) = x^{m-1} + \sum_{j=0}^{m-3} c_j t^{(m-1-j)q} x^j \tag{2.2}$$

where

$$\begin{cases} c_j = 0 & \text{if } (m-1-j)q > \tau_j \\ \beta_j(t) = c_j t^{(m-j-1)q}(1 + O(t)) & \text{if } (m-1-j)q = \tau_j. \end{cases}$$

If $q = \infty$ then $\Theta(x,t) \equiv \tilde{\Theta}(x,t) \equiv x^{m-1}$. In either case, the normalized homogeneous polynomial $Q$ associated to $L$ is

$$Q(x,z) = x^{m-1} + \sum_{j=0}^{m-3} c_j z^{m-1-j} x^j.$$

If $m = 2$ then necessarily $q = \infty$, and $L$ is analytic hypoelliptic by the well known theorems on the symplectic case.

To complete the proof of Lemma 1.3, it remains only to show that $q, Q$ are invariant in the required sense. Assume that $\{X_1, X_2\}$ is of finite type, but is linearly dependent at the origin. Let $(x,t)$, $\tilde{\Theta}$, $q$ and $\Theta$ be as defined above. Denote by $(y,s)$ some other system of coordinates, with the same origin as $(x,t)$. We regard $y, s$ both as coordinates, and as functions of $(x,t)$.

The type $m$ at $0$ is an invariant, being intrinsically defined as 1 if $X_1, X_2$ are independent at $0$, and otherwise as the minimal length of any Lie bracket of $X_1, X_2$ that is not in their span at $0$. Suppose that the span of $\{X_1, X_2\}$ equals the span of $\{\partial_y, \tilde{\Phi}(y,s)\partial_s\}$ in a neighborhood of $0$, and that $\tilde{\Phi} = y^{m-1} + \sum_{j \le m-3} \gamma_j(s) y^j$ where each $\gamma_j(0) = 0$. Then $\tilde{\Phi}(y,s) = h(y,s)\tilde{\Theta}(x,t)$ for some analytic function $h$ that does not vanish at $0$. Our task is to show that the polynomial $\Phi$ associated to $\tilde{\Phi}$ as in (2.2) is identical to $\Theta$.

Observe that $s = ct + O(t^2, xt, x^m)$ for some $c \ne 0$. Indeed, consider any monomial $D = V_1 V_2 \ldots V_n$ where each $V_i$ has smooth coefficients and belongs to the span of $X_1, X_2$ at every point of a neighborhood of $0$. To avoid confusion of notation define a function $f$ by $f(x,t) = s$. If $n < m$ then $D(f)$ must vanish at the origin, as is seen by writing each $V_i$ as a linear combination of $\partial_y$ and of $\tilde{\Phi}(y,s)\partial_s$. Choosing $D$ to equal $\partial_x^n$ in turn for each $n < m$ forces $s = f(x,t)$ to have the required form for some constant $c$. Since the differential of $f$ cannot vanish, $c$ must be nonzero.

The situation is symmetric, so we have likewise that $t = cs + O(s^2, ys, y^m)$ for some nonzero $c$. A further consequence is that $x = ay + O(s, ys, y^2)$ for some $a \ne 0$. Consider first the case where $x = ay + bs^r + O(y^2, ys, s^{r+1})$ for some $r \ge 1$ and $b \ne 0$.

Consider the subcase where $r < q$. To each monomial $y^\alpha s^\beta$ assign weight $\alpha + r^{-1}\beta$. Expand $h \cdot \tilde\Theta$ as a formal Taylor series in powers of $y, s$ at 0, and consider the weights of those monomials in this expansion that might have nonzero coefficients. For $0 \le j \le m - 3$,

$$x^{m-j-1} = (ay + bs^r + O(y^2, ys, s^{r+1}))^{m-j-1}.$$

When this is expanded, no terms of weight $< m - 1 - j$ arise, and the sum of all terms of weight equal to $m - 1 - j$ is $(ay + bs^r)^{m-j-1}$. Since $\beta_j(t) = O(t^{jq})$, each monomial in its expansion has weight $> jq/r$, and hence each monomial in the expansion of $\beta_j(t)x^{m-j-1}$ has weight $\ge (m-1-j) + jq/r$, which is strictly greater than $m - 1$. Consequently the same goes for the expansion of $h \cdot \beta_j(t)x^{m-j-1}$, and similar reasoning gives the same conclusion for $(h - h(0)) \cdot x^{m-1}$. Finally $x^{m-1} = (ay + bs^r)^{m-1}$ plus monomials of weights strictly greater than $m - 1$, and hence $h \cdot x^{m-1} = h(0)(ay + bs^r)^{m-1}$ plus monomials of weights $> m - 1$. Thus in the Taylor expansion of $\tilde\Phi = h\tilde\Theta$, no terms of weight $< m - 1$ occur, and the sum of all monomials of weight $m - 1$ is exactly $h(0)(ay + bs^r)^{m-1}$. Therefore if we expand $\tilde\Phi$ instead in powers of $y$ with coefficients depending on $s$, the coefficient of $y^{m-2}$ must equal $(m-1)a^{m-2}bh(0)s^r + O(s^{r+1})$. The numerical factor $(m-1)a^{m-2}bh(0)$ is nonzero, and therefore $\tilde\Phi$ cannot be a normalized polynomial of degree $m - 1$ in $(y, s)$, a contradiction. Thus $r$ cannot be less than $q$.

If $r = q$, the same reasoning leads to the conclusion that again, no monomials of weights $< m - 1$ occur, and that the sum of all monomials of weight $m - 1$ is $h(0)\Theta(ay + bs^q, cs)$. The coefficient of $y^{m-2}$ in this polynomial is still $(m-1)a^{m-2}bh(0)s^q$, so the same contradiction is reached.

The only remaining cases are where either $r > q$, or $x = ay + O(y^2, ys, s^N)$ for all $N$. Define the weight of any monomial $y^\alpha s^\beta$ now to be $\alpha + q^{-1}\beta$. Then the same reasoning as above leads to the conclusion that no monomials of weight $< m - 1$ occur with nonzero coefficients in the Taylor expansion of $\tilde\Phi$, and moreover that the sum of all terms of weight $m - 1$ equals $h(0)\Theta(ay, cs)$.

Writing $\tilde\Phi(y, s) = y^{m-1} + \sum_{j \le m-3} \gamma_j(s)y^j$, it follows that $\gamma_j(s) = O(s^{(m-1-j)q})$ for all $j$. Therefore if $\tilde q$ denotes the quantity obtained by applying the definition of $q$ in the coordinates $(y, s)$, we find that $\tilde q \ge q$. The situation is symmetric, so likewise $q \le \tilde q$. $\qquad\qquad\square$

## 3. Analytic preliminaries

Let $L$ be given, and adopt coordinates $(x, t)$ as described in Lemma 1.3. Let $m$ be the type at $p = 0$. Assume that $L$ is not elliptic at 0, and moreover that $q < \infty$, which implies that $m \geq 3$. To simplify notation let $E = E(Q) = NE(L, 0)$. Let $E^{1/q}$ be the set of all $q$-th roots of elements of $E(Q)$; recall that $q$ is rational.

Write

$$L = P(x, t, \partial_x, \tilde{\Theta}(x, t)\partial_t)$$

where $P$ is a $C^\omega$ function of $(x, t)$ and is an elliptic polynomial of degree 2 with respect to $\partial_x, \tilde{\Theta}\partial_t$. More concretely,

$$L = \left(a_{1,1}\partial_x + a_{1,2}\tilde{\Theta}\partial_t\right)^2 + \left(a_{2,1}\partial_x + a_{2,2}\tilde{\Theta}\partial_t\right)^2$$

for certain real valued, $C^\omega$ coefficients $a_{i,j}$, where the two by two matrix with entries $a_{i,j}(x, t)$ is nonsingular for every $(x, t)$ sufficiently close to 0.

Fix an entire holomorphic function $\Psi$ of one complex variable, not identically zero, satisfying $\Psi(z) = O(\exp(C|z|^m))$ for all $z \in \mathbb{C}$ and $\Psi(z) = O(\exp(-|z|^m))$ for all $z$ in some conic neighborhood of $\mathbb{R}$; when $m$ is even, it suffices to set $\Psi(z) = \exp(-Cz^m)$. Denote by $a, b \in \mathbb{R}$ constants to be chosen in the course of the proof of Theorem 1.4, and set $\psi(x) = e^{iax}\Psi(x - b)$. For each large $\lambda \in \mathbb{R}^+$ set

$$F_\lambda(x, t) = e^{i\lambda^m t}\psi(\lambda x)$$

**Lemma 3.1.** *Let $L$ be a locally solvable, $C^\infty$ hypoelliptic linear partial differential operator. Suppose that $L$ is analytic hypoelliptic in some neighborhood of 0. Then there exist open sets $U \subset \mathbb{R}^2$ and $\tilde{U} \subset \mathbb{C}^2$ containing 0 such that for any $N$ there exist $C, M < \infty$ such that for every $\lambda \geq 1$, there exists $G_\lambda \in C^\infty(U)$ satisfying*

$$LG_\lambda = F_\lambda$$

$$\|G_\lambda\|_{C^N(U)} \leq C\lambda^M$$

*$G_\lambda$ extends to a holomorphic function in $\tilde{U}$,*

*and*

$$|G_\lambda(x, t)| \leq C\lambda^M e^{C\lambda^m |\operatorname{Im}(x,t)|} \quad \text{for all } (x, t) \in \tilde{U}.$$

This type of result originates in work of Oleĭnik and Radkevič [16]; see the proof of Lemma 4.1 of [3] or that of Theorem 1 of [16]. □

Set

$$v_\lambda(x, t) = e^{-i\lambda^m t}G_\lambda(x, t),$$

and

$$y = \lambda x, \qquad s = \lambda^{1/q}t. \tag{3.1}$$

The final inequality in Lemma 3.1 is also satisfied by $v_\lambda$ and by each of its partial derivatives. The equation satisfied by $v_\lambda$ is

$$(e^{-i\lambda^m t} \circ L \circ e^{i\lambda^m t})v_\lambda(x,t) = \psi(\lambda x) = \psi(y).$$

**Lemma 3.2.**

$$\lambda^{m-1}\tilde{\Theta}(\lambda^{-1}y, \lambda^{-1/q}s) = \Theta(y,s) + R_\lambda(y,s)$$

*where for each $\varepsilon > 0$ there exist $\delta > 0$ and $C < \infty$ such that for all sufficiently large $\lambda \in \mathbb{R}^+$, for all $|y| \leq \lambda^{1-\varepsilon}$ and $s$ in any fixed complex neighborhood of $0$, for any $\alpha, \beta \geq 0$,*

$$|\partial_y^\alpha \partial_s^\beta R_\lambda(y,s)| \leq \begin{cases} C\lambda^{-\delta}(1+|y|)^{m-1-\alpha} & \text{if } m-1-\alpha \geq 0 \\ C\lambda^{-\delta} & \text{if } m-1-\alpha < 0. \end{cases}$$

This follows directly from the fact that $\tilde{\Theta} = \Theta$ plus terms of higher weight.  □

Define

$$B_\lambda = \lambda^{-2}e^{-i\lambda^m t} \circ L \circ e^{i\lambda^m t}.$$

Thus

$$\begin{aligned} B_\lambda &= \lambda^{-2}P\Big(\lambda^{-1}y, \lambda^{-1/q}s, \lambda\partial_y, \tilde{\Theta}(\lambda^{-1}y, \lambda^{-1/q}s)(i\lambda^m + \lambda^{1/q}\partial_s)\Big) \\ &= P\Big(\lambda^{-1}y, \lambda^{-1/q}s, \partial_y, \lambda^{-1}\tilde{\Theta}(\lambda^{-1}y, \lambda^{-1/q}s)(i\lambda^m + \lambda^{1/q}\partial_s)\Big) \\ &= P\Big(\lambda^{-1}y, \lambda^{-1/q}s, \partial_y, \lambda^{-m}(\Theta + R_\lambda)(y,s)(i\lambda^m + \lambda^{1/q}\partial_s)\Big) \\ &= P\Big(\lambda^{-1}y, \lambda^{-1/q}s, \partial_y, i(\Theta + R_\lambda)(1 - i\lambda^{-p}\partial_s)\Big) \end{aligned}$$

where

$$p = m - q^{-1}. \tag{3.2}$$

Since $q^{-1}$ equals $(m-1-j)/\tau_j$ for some $0 \leq j \leq m-3$, and since each $\tau_j \geq 1$, we have $q^{-1} \leq m-1-j \leq m-1$ and hence[3] $p \geq 1$. Setting

$$u_\lambda(y,s) = \lambda^2 v_\lambda(x,t) \tag{3.3}$$

we arrive at the equation

$$B_\lambda u_\lambda(y,s) = \psi(y). \tag{3.4}$$

Define

$$A_s = P(0, 0, \partial_y, i\Theta(y,s)). \tag{3.5}$$

---

[3] In the particular case treated in [3], $p$ was equal to 1. The main difference between the proof of Theorem 1.4 and the analysis of [3] is merely that $\lambda$ is systematically replaced by $\lambda^p$.

Our plan is essentially to analyze the equation $B_\lambda u_\lambda = \psi$ by a Neumann series argument, writing formally $B_\lambda^{-1} = \sum_0^\infty (-1)^j [A_s^{-1}(B_\lambda - A_s)]^j A_s^{-1}$. Thus we require bounds on both $A_s^{-1}$ and $B_\lambda - A_s$.

The remainder of this section summarizes properties of $A_s^{-1}$. The operator $A_s$ may be regarded in either of two ways, as an *ordinary* differential operator, acting on functions of $y$ and depending on the parameter $s \in \mathbb{C}$, or as acting on functions of both variables $(y, s)$. For the present we view it in the former light. $E^{1/q} = NE(L, 0)^{1/q}$ is thus the usual set of nonlinear eigenvalues for the family of ordinary differential operators $\zeta \mapsto A_\zeta$.

Set $\langle x \rangle = (1 + |x|^2)^{1/2}$ for any $x \in \mathbb{C}$. For $\rho \in \mathbb{R}$ and $k \in \{0, 1, 2\}$ consider the Sobolev space $\mathcal{H}_\rho^k = \mathcal{H}_\rho^k(\mathbb{R})$ of (equivalence classes of) measurable functions $f : \mathbb{R} \mapsto \mathbb{C}$ for which the following norms are finite:

$$\|f\|_{\mathcal{H}_\rho^0}^2 = \int_{\mathbb{R}} |f(x)|^2 \langle x \rangle^{-2(m-1)} e^{\rho |x|^m} dx$$

$$\|f\|_{\mathcal{H}_\rho^1}^2 = \int_{\mathbb{R}} \left( \langle x \rangle^{-2(m-1)} |\partial_x f|^2 + |f|^2 \right) e^{\rho |x|^m} dx$$

$$\|f\|_{\mathcal{H}_\rho^2}^2 = \int_{\mathbb{R}} \left( \langle x \rangle^{-2(m-1)} |\partial_x^2 f|^2 + |\partial_x f|^2 + \langle x \rangle^{2(m-1)} |f|^2 \right) e^{\rho |x|^m} dx.$$

Likewise for any open set $D \subset \mathbb{C}$ define the spaces $\mathcal{H}_\rho^k(\mathbb{R} \times D)$ to consist of all (equivalence classes of) functions of $(x, z) \in \mathbb{R} \times D$ that are holomorphic with respect to $z$, for which the following norms are finite:

$$\|f\|_{\mathcal{H}_\rho^0}^2 = \iint_{\mathbb{R} \times D} |f(x)|^2 \langle x \rangle^{-2(m-1)} e^{\rho |x|^m} dx \, dz d\bar{z}$$

$$\|f\|_{\mathcal{H}_\rho^1}^2 = \iint_{\mathbb{R} \times D} \left( \langle x \rangle^{-2(m-1)} |\partial_x f|^2 + |f|^2 \right) e^{\rho |x|^m} dx \, dz d\bar{z}$$

$$\|f\|_{\mathcal{H}_\rho^2}^2 = \iint_{\mathbb{R} \times D} \left( \langle x \rangle^{-2(m-1)} |\partial_x^2 f|^2 + |\partial_x f|^2 + \langle x \rangle^{2(m-1)} |f|^2 \right) e^{\rho |x|^m} dx \, dz d\bar{z}.$$

**Lemma 3.3.** *For any normalized homogeneous polynomial $Q$ having real coefficients, the set $E(Q)$ is discrete, and $E(Q) \cap \mathbb{R} = \emptyset = (E(Q)^{1/q}) \cap \mathbb{R}$.*

*Proof.* That $E(Q) \cap \mathbb{R} = \emptyset$ follows from positivity of the ordinary differential operator, for $\zeta \in \mathbb{R}$. Indeed, for any $f$ in the Schwartz class,

$$\int_{\mathbb{R}} (-\partial_x^2 f + Q^2(x, \zeta) f) \cdot \bar{f} = \int_{\mathbb{R}} [|f'|^2 + Q(x, \zeta)^2 |f|^2] \, dx \geq c_\zeta \|f\|_{L^2(\mathbb{R})}^2,$$

because $Q$ has real coefficients and does not vanish identically as a function of $x$ for any $\zeta$. Hence $E \cap \mathbb{R} = \emptyset$. The same reasoning applies to the set $E^{1/q}$ of nonlinear eigenvalues associated to the family of operators $A_z$, because $\Theta(x, z)$ is likewise real valued whenever $x, z$ are real.

As was proved in [5] for similar families of operators, $E(Q)$ is equal to the set of all zeros of an entire holomorphic function; see also Section 6 below. Since $E(Q) \neq \mathbb{C}$, $E(Q)$ must be discrete. $\qquad\square$

**Lemma 3.4.** *For each compact set $K \subset \mathbb{C}\backslash E^{1/q}$ there exists $r > 0$ such that for every $\rho \in [-r, r]$ and $\zeta \in K$, $A_\zeta : \mathcal{H}_\rho^2 \mapsto \mathcal{H}_\rho^0$ is an isomorphism, whose inverse is bounded uniformly in $\zeta \in K$, $|\rho| \leq r$.*

*Proof.* This follows from the method of [3], once it is shown that the (two dimensional) nullspace of $A_\zeta$ contains no functions that are $O(\exp(-\delta|x|^m))$ for some $\delta > 0$, since any function in the nullspace of $A_\zeta$ either decays at such a rate as $|x| \to \infty$, or tends to $\infty$ in modulus as $x \to +\infty$ or as $x \to -\infty$. The hypothesis that $K$ does not intersect $E^{1/q}$ means that $\partial_x^2 - \Theta(x, \zeta)^2$ has no such solutions, for any $\zeta \in K$.

To relate $A_\zeta$ to this operator, write $A_\zeta = \alpha \partial_x^2 + i\beta(\partial_x \circ \Theta(x, \zeta) + \Theta(x, \zeta)\partial_x) - \gamma\Theta(x, \zeta)^2$ where $\alpha, \beta, \gamma \in \mathbb{R}$, $\alpha, \gamma$ are positive, $\beta^2 < \alpha\gamma$, and $\Theta$ denotes both a function and the operator defined by multiplication by that function. Fix a polynomial $R(x, z)$, with real coefficients, satisfying $\partial R/\partial x = Q(x, z)$, and consider

$$e^{i\beta R(x,\zeta)/\alpha} \circ A_\zeta \circ e^{-i\beta R(x,\zeta)/\alpha} = \alpha\partial_x^2 + (\alpha^{-1}\beta^2 - \gamma)\Theta(x, \zeta)^2.$$

The coefficient $\alpha^{-1}\beta^2 - \gamma$ is negative. Dilating both variables $x, \zeta$ by appropriate factors and multiplying the operator by a constant reduces it to $\partial_x^2 - \Theta(x, \zeta)^2$.

Now the factor $\exp(i\beta R(x, \zeta)/\alpha)$ may not be a bounded function of $x$ when $\zeta$ is not real, but since $R$ has real coefficients and $\beta/\alpha$ is real, this factor and its inverse are $O(\exp(C|x|^{m-2}))$ for $\zeta$ in any compact set and $x \in \mathbb{R}$. Therefore since $\partial_x^2 - \Theta(x, \zeta)^2$ has no solutions bounded by $\exp(-\delta|x|^m)$ as $x \to \pm\infty$, no matter how small the exponent $\delta$, neither has $A_\zeta$. $\qquad\square$

The proofs of the remaining lemmas of this section are essentially identical to those in Section 3 of [3], and will not be repeated here. Define $\Psi_{a,b}(x) = \exp(iax)\Psi(x - b)$.

**Lemma 3.5.** *Suppose that $\Psi$ is a Schwartz function that does not vanish identically. For each $\zeta_0 \in E^{1/q}$, for any circle $\Gamma$ centered at $z_0$ whose closure contains no other points of $E^{1/q}$, there exist $\sigma \in \{0, 1, 2, \dots\}$, $a, b \in \mathbb{R}$ and $\varphi \in C_0^\infty(\mathbb{R})$ such that*

$$\int_{\mathbb{R}} \oint_\Gamma \varphi(x) A_\zeta^{-1} \Psi_{a,b}(x) \, \zeta^\sigma d\zeta \, dx \neq 0.$$

The proof relies on the fact that the span of all functions $\Psi_{a,b}$ is dense in $L^2(\mathbb{R})$. Otherwise it is essentially the same as the proof of Lemma 3.2 of [3].

For $\zeta, \tau \in \mathbb{C}$ define

$$A_{\zeta,\tau} = P(0, 0, \partial_x, i\tau\Theta(x, \zeta))$$

Then for $0 \neq \tau$ in a small conic neighborhood of $\mathbb{R}$,

$$A_{\zeta,\tau} = P(0, 0, \partial_x, \tau^{1/m} i\Theta(\tau^{1/m}x, \tau^{1/mq}\zeta))$$
$$= \tau^{2/m} P(0, 0, \partial_w, i\Theta(w, \tau^{1/mq}\zeta))$$

where $w = \tau^{1/m}x$.

**Lemma 3.6.** *If $\tau$ belongs to a sufficiently small conic neighborhood of $\mathbb{R}$ then for any $\zeta \in \mathbb{C}$, $A_{\zeta,\tau}$ annihilates some nonzero Schwartz class function if and only if $\tau^{1/mq}\zeta \in E^{1/q}$. If $\tau^{1/mq}\zeta \notin E^{1/q}$ then $A_{\zeta,\tau} : \mathcal{H}^2_\rho \mapsto \mathcal{H}^0_\rho$ is invertible, for all sufficiently small $|\rho|$. This holds uniformly for any compact set of such $(\zeta, \tau)$.*

**Lemma 3.7.** *There exists a conic neighborhood of $\mathbb{R}$ that is disjoint from $E^{1/q}$.*

## 4. Necessity: Reduction to eigenvalue problems

This section contains the core of the proof of Theorem 1.4, which asserts that analytic hypoellipticity fails to hold whenever a set of eigenvalues is nonempty. Recall that $B_\lambda u_\lambda = \psi$, and

$$|u_\lambda| + |\nabla u_\lambda| + |\nabla^2 u_\lambda| \leq \lambda^M \exp(C\lambda^m |\operatorname{Im}(t)|) \leq \exp(C_0\lambda^p) \tag{4.1}$$

for all $|y| \leq c\lambda$ and $|s| \leq C_1$, for some $C_0 < \infty$, where $(y, s) \in \mathbb{R} \times \mathbb{C}$ and $C_1$ may be taken to be as large as we wish provided that $\lambda$ and $C_0$ are sufficiently large. Here $\nabla$ denotes the gradient in both variables $y, s$.

Assume that $E^{1/q}$ is nonempty, as hypothesized in Theorem 1.4, and fix any $\zeta_0 \in \bar{E}^{1/q}$ having strictly negative imaginary part. $E = \bar{E}$ because the coefficients of $Q$ are real, and $E^{1/q} \cap \mathbb{R} = \emptyset$, so such a point exists. Define $\Gamma = \{|s - \zeta_0| = r\}$ to be a circle centered at $\zeta_0$, contained in the open lower half plane, such that the intersection of $E^{1/q}$ with the closed disk bounded by $\Gamma$ contains only the point $\zeta_0$.

Consider $(y, s) \in \mathbb{R} \times \Gamma$, and write $s = \zeta_0 + re^{i\theta}$. When acting on restrictions to $\mathbb{R} \times \Gamma$ of functions holomorphic with respect to $s$, $B_\lambda$ takes the form

$$B_\lambda = P\Big(\lambda^{-1}y, \lambda^{-1/q}(\zeta_0 + re^{i\theta}), \partial_y, i[\Theta + R_\lambda](1 - i\lambda^{-p}(-ir^{-1}e^{-i\theta}\partial_\theta))\Big)$$

where $\Theta + R_\lambda$ is evaluated at $(y, \zeta_0 + re^{i\theta}) = (y, s)$. Denote by $B_\lambda^*$ the transpose of $B_\lambda$ in $L^2(\mathbb{R} \times \Gamma, dyd\theta)$, where functions of $s \in \Gamma$ are identified with periodic functions of $\theta$, $|\theta| \leq \pi$. Thus

$$B_\lambda^* = P^*\Big(\lambda^{-1}y, \lambda^{-1/q}s, -\partial_y, i(1 - i\lambda^{-p}(ir^{-1}\partial_\theta \circ e^{-i\theta})) \circ [\Theta + R_\lambda]\Big)$$

where

$$P^*(x, t, -\partial_x, -\partial_t \circ \tilde{\Theta}) = \Big(-\partial_x \circ a_{1,1} - \partial_t \circ \tilde{\Theta}a_{1,2}\Big)^2 + \Big(-\partial_x \circ a_{2,1} - \partial_t \circ \tilde{\Theta}a_{2,2}\Big)^2,$$

and where $a_{i,j}$ denotes also the operator defined by multiplication by the function $a_{i,j}$.

Let $\sigma, \varphi$ be as in Lemma 3.5. The proofs of the next two lemmas will be discussed at the end of this section.

**Lemma 4.1.** *There exist $\rho, \delta, \varepsilon, C \in \mathbb{R}^+$ such that for each sufficiently large $\lambda \in \mathbb{R}^+$ there exists $f : \mathbb{R} \times \Gamma \mapsto \mathbb{C}$, supported where $|y| \leq \lambda^{1-\varepsilon}$, satisfying*

$$\iint_{\mathbb{R} \times \Gamma} |B_\lambda^* f - s^\sigma (\partial s / \partial \theta) \varphi(y)|^2 e^{\rho |y|^m} \, dy \, d\theta \qquad \leq e^{-\delta \lambda^p} \qquad (4.2)$$

$$\iint_{\mathbb{R} \times \Gamma} \left( |f|^2 + |\nabla f|^2 \right) e^{\rho |y|^m} \, dy \, d\theta \qquad \leq C \qquad (4.3)$$

$$\iint_{\mathbb{R} \times \Gamma} |f(y, \theta) - s^\sigma (\partial s / \partial \theta) A_s^{-1} \varphi(y)|^2 e^{\rho |y|^m} \, dy \, d\theta \leq C \lambda^{-\delta} \qquad (4.4)$$

*where $s = \zeta_0 + r e^{i\theta}$.*

The next lemma is conditional in nature; we will soon see that $u_\lambda$ cannot satisfy its conclusion, hence must not satisfy its hypotheses.

**Lemma 4.2.** *Suppose that $B_\lambda u_\lambda = \psi$ and that $u_\lambda$ satisfies the bound (4.1). Then for any $C_1 < \infty$ and $\mu, \varepsilon, \varrho > 0$ there exists $C < \infty$ such that for all sufficiently large $\lambda$ and all $|s| \leq C_1$,*

$$\int_{|y| \leq \lambda^{1-\varepsilon}} \left[ |u_\lambda(y, s)|^2 + |\nabla u_\lambda(y, s)|^2 \right] e^{-\varrho |y|^m} \, dy \leq e^{\mu \lambda^p}. \qquad (4.5)$$

Theorem 1.3 may now be proved, as follows. Suppose that $L$ were analytic hypoelliptic in some neighborhood of the origin. Define $F_\lambda$ as in Section 3, solve $LG_\lambda = F_\lambda$ as in Lemma 3.1, and define first $u_\lambda$ and then $v_\lambda$ as was done following Lemmas 3.1 and 3.2, respectively. Then Lemma 4.2 yields a strong bound for $u_\lambda$. Let $f$ be an approximate solution for the adjoint operator $B_\lambda^*$, satisfying the conclusions of Lemma 4.1. Define

$$\omega = \int_{\mathbb{R}} \oint_\Gamma \varphi(y) u_\lambda(y, s) s^\sigma \, ds \, dy.$$

Then $\omega = 0$ since $u_\lambda$ is holomorphic with respect to $s$. We will prove that

$$\omega = c_0 + O(\lambda^{-\delta}) \qquad \text{as } \lambda \to +\infty \qquad (4.6)$$

for some $c_0 \neq 0$ and $\delta > 0$, thus arriving at a contradiction.

Choosing $\mu$ in Lemma 4.2 to be less than $\delta$ in Lemma 4.1, and $\varrho$ in Lemma 4.2 to be less than $\rho$ in Lemma 4.1,

$$\omega = \iint_{\mathbb{R} \times \Gamma} [\varphi(y) s^\sigma (\partial s / \partial \theta)] u_\lambda(y, s) \, d\theta \, dy$$
$$= \iint_{\mathbb{R} \times \Gamma} (B_\lambda^* f)(y, s) u_\lambda(y, s) \, d\theta \, dy + O(e^{-\varepsilon \lambda^p})$$

for some $\varepsilon > 0$, by (4.2) and (4.5). Because $f$ has compact support with respect to $y$ and $\Gamma$ has no boundary, it is permissible to integrate by parts to rewrite the last line as

$$\omega = \iint_{\mathbb{R} \times \Gamma} f(y, s) B_\lambda u_\lambda(y, s) \, d\theta \, dy + O(e^{-\varepsilon \lambda^p}).$$

Note that for each $s \in \mathbb{C}$, $A_s$ is its own transpose. Therefore

$$\omega = \iint_{\mathbb{R} \times \Gamma} f(y, s) \psi(y) \, d\theta \, dy + O(e^{-\varepsilon \lambda^p})$$
$$= \iint_{\mathbb{R} \times \Gamma} \left[ s^\sigma (\partial s / \partial \theta) A_s^{-1} \varphi(y) \right] \psi(y) \, d\theta \, dy + O(\lambda^{-\delta})$$
$$= \oint_\Gamma \left( \int_{\mathbb{R}} A_s^{-1} \varphi(y) \cdot (A_s A_s^{-1} \psi)(y) \, dy \right) s^\sigma \, ds + O(\lambda^{-\delta})$$
$$= \oint_\Gamma \int_{\mathbb{R}} \varphi(y) A_s^{-1} \psi(y) \, dy s^\sigma \, ds + O(\lambda^{-\delta})$$
$$= c_0 + O(\lambda^{-\delta})$$

where $c_0 \neq 0$ by Lemma 3.5.                                                     $\square$

We next prove Lemma 4.1. It is convenient for $B_\lambda^*$ to be globally defined, on the Cartesian product of $\mathbb{R}$ with a complex neighborhood of $\Gamma$. This can be accomplished as in the first paragraph of the proof of Lemma 5.1 of [3]: Fix a constant $\varepsilon_0 > 0$ satisfying

$$(1 - \varepsilon_0) m > p = m - q^{-1}. \tag{4.7}$$

Extend the coefficients $a_{i,j}(\lambda^{-1} y, \lambda^{-1/q} s)$ to $\tilde{a}_{i,j}(\lambda^{-1} y, \lambda^{-1/q} s)$ so that $\tilde{a}_{i,j}(\lambda^{-1} y, \lambda^{-1/q} s)$ is identically equal to $a_{i,j}(\lambda^{-1} y, \lambda^{-1/q} s)$ when $|y| \leq \lambda^{1-\varepsilon_0}$ and is equal to $a_{i,j}(0, \lambda^{-1/q} s)$ when $|y| \geq 2\lambda^{1-\varepsilon_0}$, and so that for $\lambda^{1-\varepsilon_0} \leq |y| \leq 2\lambda^{1-\varepsilon_0}$,

$$\tilde{a}_{i,j}(\lambda^{-1} y, \lambda^{-1/q} s) - a_{i,j}(\lambda^{-1} y, \lambda^{-1/q} s) = O(\lambda^{-1} y) = O(\lambda^{-\varepsilon_0})$$

and

$$\partial_s^\beta \partial_y^\alpha \tilde{a}_{i,j} = O(\lambda^{-\alpha - \beta/q}) \text{ for all } \alpha, \beta \geq 0.$$

Similarly, $R_\lambda$ may be extended so that

$$|\partial_y^\alpha \partial_s^\beta R_\lambda| \leq C \lambda^{-\varepsilon_0} (1 + |y|)^{\max(0, m-1-\alpha)}$$

for all $|y| \leq \lambda^{1-\varepsilon_0}$, $|s| \leq C$.

Let $\Lambda \in \mathbb{R}^+$ be a large constant to be chosen below, and given any large $\lambda$, construct concentric open annuli $\Omega_0 \supset \Omega_1 \supset \cdots \supset \Omega_{N+1} = \Omega \supset \Gamma$ so that

1. $E^{1/q} \cap \overline{\Omega_0} = \emptyset$,
2. The three sets $\Omega_0, \Omega_1, \Omega$ are all independent of $\lambda$,
3. $N \geq c\Lambda^{-1}\lambda^p$, and
4. distance $(\partial\Omega_j, \Omega_{j+1}) \geq \Lambda\lambda^{-p}$ for all $j$.

Set $\phi(y,s) = s^\sigma(\partial s/\partial\theta)\varphi(y) = is^\sigma(s - \zeta_0)\varphi(y)$.

Decompose

$$B_\lambda^* = A_s + K + \mathcal{E}$$

where

$$A_s = P(0, 0, \partial_y, i\Theta(y,s)),$$

$$A_s + K = P^*(\lambda^{-1}y, \lambda^{-1/q}s, -\partial_y, i\Theta + iR_\lambda).$$

Thus the remainder $\mathcal{E}$ represents all contributions of $\partial_s$. $A_s$ and $K$ will be regarded both as ordinary differential operators depending holomorphically on the parameter $s$, and as operators acting on functions of $(y,s)$ that are holomorphic with respect to $s$.

For all sufficiently small $\rho > 0$, $A_s : \mathcal{H}_\rho^2(\mathbb{R}) \mapsto \mathcal{H}_\rho^0(\mathbb{R})$ is an isomorphism, uniformly for all $s \in \overline{\Omega_0}$, by Lemma 3.4. Since $A_s^{-1}$ depends holomorphically on $s \in \mathbb{C}\backslash E^{1/q}$, it follows that $A_s : \mathcal{H}_\rho^2(\mathbb{R} \times U) \mapsto \mathcal{H}_\rho^0(\mathbb{R} \times U)$ is also an isomorphism for every open set $U \subset \Omega_0$, uniformly in $\rho, U$ provided that $\rho$ is sufficiently small.

By Cauchy's inequality, $\partial_s : \mathcal{H}_\rho^k(\mathbb{R} \times \Omega_j) \mapsto \mathcal{H}_\rho^k(\mathbb{R} \times \Omega_{j+1})$ has norm $O(\Lambda^{-1}\lambda^p)$, uniformly in all parameters, because the distance between the boundaries is at least $\Lambda\lambda^{-p}$, so

$$\lambda^{-p}\partial_s : \mathcal{H}_\rho^k(\mathbb{R} \times \Omega_j) \mapsto \mathcal{H}_\rho^k(\mathbb{R} \times \Omega_{j+1}) \text{ has norm } \leq C\Lambda^{-1} \text{ for all } j, \lambda.$$

Moreover, the norm is $O(\lambda^{-p})$ for $j = 0$, since the distance from $\partial\Omega_0$ to $\partial\Omega_1$ is independent of $\lambda$. From this, from the definitions of the norms, and from straightforward computation there follows the first conclusion of the following lemma.

**Lemma 4.3.** *The remainder terms $K, \mathcal{E}$ satisfy the following bounds.*

$$\mathcal{E} : \mathcal{H}_\rho^2(\mathbb{R} \times \Omega_j) \mapsto \mathcal{H}_\rho^0(\mathbb{R} \times \Omega_{j+1}) \text{ is } O(\Lambda^{-1} + \lambda^{-\varepsilon})$$

*for some $\varepsilon > 0$, as $\lambda \to +\infty$. For $j = 0$, it is $O(\lambda^{-\varepsilon})$. Secondly,*

$$K : \mathcal{H}_\rho^2(\mathbb{R} \times \Omega_j) \mapsto \mathcal{H}_\rho^0(\mathbb{R} \times \Omega_j) \text{ is } O(\lambda^{-\varepsilon})$$

*for some $\varepsilon > 0$, as $\lambda \to +\infty$. These bounds hold for all $|\rho| \leq \rho_0 \ll 1$, uniformly in $\lambda, j$.*

The second conclusion follows from the bounds of Lemma 3.2 on $R_\lambda$ and its derivatives, and the presence of negative powers of $\lambda$ in the expressions $\lambda^{-1}y, \lambda^{-1/q}s$.

Fix $\rho > 0$ sufficiently small that Lemma 3.4 applies for all $\zeta \in \overline{\Omega_0}$. As an approximate solution to $B_\lambda^* f \approx \phi$ define

$$f = \sum_{j=0}^{N} (-1)^j [A_s^{-1}(K + \mathcal{E})]^j A_s^{-1} \phi.$$

Then

$$B_\lambda^* f - \phi = \pm[(K + \mathcal{E})A_s^{-1}]^{N+1}\phi.$$

Lemma 4.3 implies

$$\|[(K + \mathcal{E})A_s^{-1}]^j \phi\|_{\mathcal{H}_\rho^0(\mathbb{R} \times \Omega_j)} \leq C^j \lambda^{-\varepsilon} \Lambda^{1-j} \tag{4.8}$$

for all $j \geq 1$, and

$$\|[A_s^{-1}(K + \mathcal{E})]^j A_s^{-1} \phi\|_{\mathcal{H}_\rho^2(\mathbb{R} \times \Omega_j)} \leq C^j \lambda^{-\varepsilon} \Lambda^{1-j}. \tag{4.9}$$

Thus if $\Lambda$ is chosen to be greater than $2C$ for the constants $C$ in (4.8) and (4.9), the fact that $N \geq c\Lambda^{-1}\lambda^p$ may be used to conclude that

$$\|B_\lambda^* f - \phi\|_{\mathcal{H}_\rho^0(\mathbb{R} \times \Omega)} \leq C(C/\Lambda)^N \leq Ce^{-\delta\lambda^p}$$

for some $\delta > 0$. (4.3) follows by summing over $j$, if $\Lambda$ is chosen to be sufficiently large. And

$$f - A_s^{-1}\phi = \sum_{j=1}^{N} (-1)^j [A_s^{-1}(K + \mathcal{E})]^j A_s^{-1} \phi = O(\lambda^{-\delta})$$

in $\mathcal{H}_\rho^2(\mathbb{R} \times \Omega)$ norm. Specializing to $s \in \Gamma$ and invoking Cauchy's inequality to majorize the $L^2(\Gamma)$ norm by that of $L^2(\Omega)$, this is (4.4).

This $f$ cannot be expected to have compact support. So fix a cutoff function $\eta \in C_0^\infty[-1, 1]$ that is $\equiv 1$ for $|y| \leq 1/2$, and define $\tilde{f}(y, s) = \eta(\lambda^{-1+\varepsilon_0}y)f(y, s)$. Then $\tilde{f}$ has all the required properties if $\rho$ is replaced by $\rho/2$ and $\varepsilon_0$ is sufficiently small. In particular, $\tilde{f}$ is supported where $|y| \leq \lambda^{1-\varepsilon_0}$, in which region the coefficients of the modified operator $B_\lambda^*$ agree with those of the original operator $B_\lambda$. In each of the three conclusions of Lemma 4.1, the contribution of $\tilde{f} - f$, in the $L^2$ norm with respect to the weight $\exp(\rho|y|^m/2)$, is $O(\exp(-c\lambda^{(1-\varepsilon_0)m}))$ by the bound (4.4) for $f, \nabla f$ with respect to the weight $\exp(\rho|y|^m)$. Since $(1 - \varepsilon_0)m > p$, this last bound is $\ll \exp(-\lambda^p)$ for large $\lambda$. $\qquad \square$

The proof of Theorem 1.4 is now complete modulo the proof of Lemma 4.2, which is outlined in Section 9.

## 5. Nonlinear eigenvalue problems

Assume that $Q(x,z) = x^{m-1} + \sum_{j=0}^{m-3} c_j z^{m-1-j} x^j$ where each $c_j \in \mathbb{R}$, $x \in \mathbb{R}$, $z \in \mathbb{C}$, and $m \geq 2$. Let $\mathcal{L}_z = -\partial_x^2 + Q(x,z)^2$, and recall that $E(Q)$ denotes the set of all nonlinear eigenvalues of the family of operators $\mathcal{L}_z$. This section contains various results connected with Conjecture 1.5, which says that $E(Q) = \emptyset$ if and only if $Q(x,z) \equiv x^{m-1}$.

**Lemma 5.1.** *If $Q(x,z) \equiv x^{m-1}$ then $E(Q) = \emptyset$.*

This is implied by other results, for $L = \partial_t^2 + x^{2(m-1)} \partial_t^2$ is analytic hypoelliptic by Theorem 1.2, which by Theorem 1.4 precludes the existence of eigenvalues.

*Proof.* $\mathcal{L}_z$ does not depend on $z$, and we have already remarked that $E(Q) \cap \mathbb{R}$ is always empty. $\qquad\square$

Before discussing the proofs of Theorem 1.6 and Proposition 1.7, we introduce a closely related class of nonlinear eigenvalue problems for which we are able to prove the analogue of Conjecture 1.5. Let $Q$ be any normalized homogeneous polynomial, with real coefficients, of odd degree $m-1$. Define

$$\mathbb{L}_z = (\partial_x + Q(x,z)) \circ (-\partial_x + Q(x,z)),$$

$$E'(Q) = \{z \in \mathbb{C} : \text{ there exists } 0 \neq f \in \mathcal{S} \text{ such that } \mathbb{L}_z f = 0\}.$$

**Theorem 5.2.** *Let $Q$ be a normalized homogeneous polynomial of odd degree $m-1$, with real coefficients. Then $E'(Q) = \emptyset$ if and only if $Q(x,z) \equiv x^{m-1}$.*

*Proof.* Introduce

$$P(x,z) = \int_0^x Q(y,z)\,dy.$$

The function

$$\psi_z^-(x) = e^{P(x,z)} \int_{-\infty}^x e^{-2P(y,z)}\,dy$$

satisfies $\mathbb{L}_z \psi_z^- \equiv 0$, and $\psi_z^-$ tends rapidly to 0 as $x \to -\infty$. On the other hand, $\exp(P(x,z))$ is a solution of $\mathbb{L}_z$ whose modulus tends rapidly to $+\infty$ as $x \to -\infty$. Therefore since $\mathbb{L}_z$ has only a two-dimensional nullspace, each solution of $\mathbb{L}_z$ either tends rapidly to $\infty$ in modulus as $x \to -\infty$, or is a scalar multiple of $\psi_z^-$. Thus $z \in E'(Q)$ if and only if $\psi_z^-(x)$ behaves as a Schwartz function as $x \to +\infty$. Defining

$$W(z) = \int_{\mathbb{R}} \exp(-2P(x,z))\,dx,$$

one has

$$E'(Q) = \{z \in \mathbb{C} : W(z) = 0\}. \tag{5.1}$$

Indeed, if $W(z) \neq 0$, then $\psi_z^-(x) \sim W(z) \exp(P(x, z))$ tends to $\infty$ in modulus as $x \to +\infty$. If $W(z) = 0$ then $\psi_z^-(x) = -\exp(P(x, z)) \int_x^\infty \exp(-2P(y, z)) dy$, which tends to 0 as $x \to +\infty$. If $Q(x, z) = (x - cz)^{m-1}$ then a change of the contour of integration reveals that $W(z) = c_1 \exp(c_2 z^m)$ for some constants $c_j$, so $W$ has no zeros.

Since $P$ is a polynomial of even degree with positive leading coefficient independent of $z$, $W$ is clearly an entire holomorphic function. For any $z \in \mathbb{C}$,

$$\begin{aligned}
|W(z)| &\leq \int_{\mathbb{R}} e^{-2x^m/m} e^{C(|z|^{m-1}|x| + |z||x|^{m-1})} dx \\
&\leq \int_{\mathbb{R}} e^{-x^m/m} e^{C|z|^m} dx \\
&= C' e^{C|z|^m},
\end{aligned}$$

so $W$ is of finite order, and its order does not exceed $m$. Any entire function of finite order with no zeros must be of the form $\exp(R(z))$ for some polynomial $R$, so in order to prove that $W$ must have zeros, it suffices merely to obtain sufficient information on the asymptotic behavior of $W$ as $\mathbb{R} \ni z \to +\infty$ to rule out $W = e^R$.

Restrict attention for the remainder of the proof to the case where $z \in \mathbb{R}^+$. Substituting $x = zy$ yields

$$W(z) = z \int_{\mathbb{R}} e^{-2z^m P(y, 1)} dy.$$

Define

$$\eta = \max_{x : Q(x, 1) = 0} -P(x, 1) = \max_{y \in \mathbb{R}} -P(y, 1).$$

Let $x_1, \ldots x_N$ be those real points at which $Q(x_j, 1) = 0$ and $-P(x_j, 1) = \eta$. Since $P$ has even degree, real coefficients, and negative leading coefficient, $-P$ has at least one global maximum, so there exists at least one such point $x_j$.

These points determine the asymptotic behavior of $W(z)$ as $z \to +\infty$. To compute the asymptotics choose $0 < \delta < \frac{1}{2} \min_{i,j} |x_i - x_j|$. Then for some $\varepsilon > 0$,

$$\int_{\mathbb{R} \setminus \cup_j [x_j - \delta, x_j + \delta]} e^{-2z^m P(y, 1)} dy = O(e^{2(\eta - \varepsilon) z^m})$$

by straightforward majorization.

For each $j$, let $k_j$ be the smallest integer for which $\partial^{k_j} P(y, 1) / \partial y^{k_j}$, evaluated at $y = x_j$, is nonzero. Since $-P(y, 1)$ has a local maximum at $x_j$, each $k_j$ must be even, and the derivative $P^{(k_j)}(x_j, 1)$ is strictly positive. By the method of real stationary phase,

$$\int_{x_j-\delta}^{x_j+\delta} e^{-2z^m P(y,1)}\, dy = e^{-2\eta z^m}\Big(c(k_j)(z^m P^{(k_j)}(x_j,1))^{-1/k_j} + O(z^{-2m/k_j})\Big),$$

where each $c(k_j)$ is strictly positive. Defining $k = \max k_j$, we find that as $\mathbb{R}^+ \ni z \to +\infty$,

$$W(z) = c e^{-2\eta z^m} z^{1-m/k}(1 + O(z^{-\varepsilon})) \tag{5.2}$$

for some $c, \varepsilon > 0$. Thus the logarithm of $W$, restricted to $\mathbb{R}^+$, cannot be a polynomial unless $1 - m/k = 0$. In that event $k = m$, that is, $y \mapsto P(y, 1)$ has a zero of order $m$ at some $x_j$. Since $P$ is a polynomial of degree $m$, this is only possible if $P(y, 1)$ takes the form $c(y - x_j)^m$. Thus because of the normalizations already imposed on $Q$, $Q(y, 1) \equiv y^{m-1}$. $\qquad\square$

Our next result concerns the ordinary differential operators $\mathcal{L}_z = -\partial_x^2 + Q(x, z)^2$, for generic $Q$. Generalize the setting by permitting the coefficients of $Q$ to be complex, and the coefficient of $x^{m-2}$ to be nonzero.[4] The definitions of $\mathcal{L}_z$ and of $E(Q)$ still make sense.

**Proposition 5.3.** *Suppose that for some $n \geq 1$, for each $0 \leq j \leq m - 2$ there is given an entire holomorphic function $g_j : \mathbb{C}^n \mapsto \mathbb{C}$. Set*

$$Q_\zeta(x, z) = x^{m-1} + \sum_{j=0}^{m-2} g_j(\zeta) z^{m-1-j} x^j,$$

$$\mathcal{L}_z^\zeta = -\partial_x^2 + Q_\zeta(x, z)^2,$$

*and*

$$E_\zeta = \{z \in \mathbb{C} : \text{ there exists } 0 \neq f \in \mathcal{S}(\mathbb{R}) \text{ such that } \mathcal{L}_z^\zeta f = 0\}.$$

*Suppose that there exists $\zeta_0 \in \mathbb{C}^n$ for which $E_{\zeta_0} \neq \emptyset$. Then*

$$E_\zeta \neq \emptyset \text{ for generic } \zeta,$$

---

[4]A normalization was carried out in Section 2 to reduce the study of analytic hypoellipticity to the case $c_{m-2} = 0$. A term $c_{m-2}zx^{m-2}$ can likewise be eliminated in the eigenvalue theory, as follows: Suppose that $Q(x, z) = x^{m-1} + \sum_{j=0}^{m-2} c_j z^{m-1-j} x^j$ where each $c_j \in \mathbb{R}$, $x \in \mathbb{R}$, $z \in \mathbb{C}$, and $m \geq 2$. Define $\tilde{Q}(x, z) = Q(x - (m-1)^{-1}c_{m-2}z, z)$, and let $\tilde{\mathcal{L}}_z$ be the associated family of Schrödinger operators. If $\mathcal{L}_\zeta f = 0$ then $f_\zeta$ extends to an entire holomorphic function of $x \in \mathbb{C}$, and if moreover $f_\zeta(x) \to 0$ as $\mathbb{R} \ni x \to \pm\infty$, then $\mathbb{R} \ni x \mapsto f(x + b)$ defines a Schwartz function, for any $b \in \mathbb{C}$ [10]. Thus $\tilde{f}(x) = f_\zeta(x - (m-1)^{-1}c_{m-2}\zeta)$ is a Schwartz function annihilated by $\hat{\mathcal{L}}_\zeta$. Therefore there is a one to one correspondence between $E(Q)$ and $E(\tilde{Q})$.

*in the sense that the set of $\zeta \in \mathbb{C}^n$ for which $E_\zeta = \emptyset$ is pluripolar.*[5]

For the normalized homogeneous polynomial $Q(x,z) = x^{m-1} - z^{m-1}$, $E(Q)$ is already known [8] to be nonempty. Therefore, taking $n = m-2$ and defining $g_{m-2} \equiv 0$ and $g_j(\zeta) = \zeta_j$ for all $0 \leq j \leq m-3$ where $\zeta = (\zeta_0, \ldots \zeta_{m-3}) \in \mathbb{C}^{m-2}$, we obtain the following corollary.

**Proposition 5.4.** *The set of all $c = (c_0, \ldots c_{m-3}) \in \mathbb{R}^{m-2}$ for which $E(Q) = \emptyset$ for $Q = x^{m-1} + \sum_{j \leq m-3} c_j z^{m-1-j} x^j$ is pluripolar.*

Combining Proposition 5.4 with Theorem 1.4, we deduce Theorem 1.6 as well. □

We turn to the detailed analysis of the nonlinear eigenvalue problem for $\mathcal{L}_z = -\partial_x^2 + Q(x,z)^2$, assuming always that $Q$ is normalized, homogeneous, and has real coefficients. Define $P(x,z) = \int_0^x Q(y,z)\, dy$.

**Lemma 5.5.** *Suppose that $m$ is even. For each $z \in \mathbb{C}$ there exist unique entire holomorphic functions $x \mapsto \psi_z^\pm(x)$, depending holomorphically on $z$ and on the coefficients $c_j$ of $Q$ for fixed $m$, satisfying*

$$\mathcal{L}_z \psi_z^\pm = 0$$

*and*

$$\begin{cases} \psi_z^+(x) & = e^{-P(x,z)}(\ \ x)^{-(m-1)/2}(1 + O(|x|^{-1})) \text{ as } \mathbb{R} \ni x \to +\infty \\ \psi_z^-(x) & = e^{-P(x,z)}(-x)^{-(m-1)/2}(1 + O(|x|^{-1})) \text{ as } \mathbb{R} \ni x \to -\infty. \end{cases} \tag{5.3}$$

*Any solution $f$ of $\mathcal{L}_z f = 0$ on the real axis either tends to $\infty$ in modulus as $x \to +\infty$, or is a scalar multiple of $\psi_z^+$. Likewise any solution either tends to $\infty$ in modulus as $x \to -\infty$, or is a scalar multiple of $\psi_z^-$. For any $Q$, the asymptotics (5.3) remain valid for all $x$ in some conic neighborhood of $\mathbb{R}$.*

*If $m$ is odd then all assertions remain valid if the asymptotics for $\psi_z^-$ are changed to*

$$\psi_z^-(x) = e^{+P(x,z)}(-x)^{-(m-1)/2}(1 + O(|x|^{-1})) \text{ as } \mathbb{R} \ni x \to -\infty.$$

*For each $Q$ there exists $C < \infty$ such that*

$$|\psi_z^+(x)| + |\partial_x \psi_z^+(x)| \leq C \exp(C|z|^m) \text{ for all } x \geq 0$$

*and*

$$|\psi_z^-(x)| + |\partial_x \psi_z^-(x)| \leq C \exp(C|z|^m) \text{ for all } x \leq 0.$$

---

[5]I am indebted for this device to a paper of D. Barrett [2]; it is a general idea in spectral theory [1].

**Lemma 5.6.** *Let $Q$ be fixed. There exists $C < \infty$ such that for all $z \in \mathbb{C}$ and all $x \geq C + C|z|$, for $k = 0, 1$,*

$$\left| \partial_x^k \psi_z^+(x) - (-1)^k e^{-P(x,z)} x^{-(m-1)/2} Q(x,z)^k \right| \leq C \left| e^{-P(x,z)} |x|^{-(m-3)/2} Q(x,z)^k \right|.$$

*The corresponding bound holds for $\psi_z^-(x)$ for all $x \leq -C - C|z|$.*

These two lemmas are proved as in [5] and as in the proofs of Lemmas 6.1 and 6.4 below, by rewriting the ordinary differential equation as a first order system, diagonalizing the system modulo a coefficient matrix that is suitably small for large $|x|$, and solving an integral equation. Holomorphic dependence on the coefficients $c_j$ is not addressed in Lemmas 6.1 and 6.4, but follows directly from their proofs. □

Central to our analysis is the Wronskian

$$W(z) = \operatorname{Det} \begin{pmatrix} \psi_z^+ & \psi_z^- \\ \partial_x \psi_z^+ & \partial_x \psi_z^- \end{pmatrix}(0). \tag{5.4}$$

The function $W$ in the proof of Theorem 5.2 may be interpreted as such a Wronskian.

**Corollary 5.7.** *$W$ has the following properties.*

1. *$W$ depends holomorphically on $z$ and on the coefficients $c_j$ of $Q$.*
2. *$z \in E(Q)$ if and only if $W(z) = 0$.*
3. *For each $Q$ there exists a constant $C$ such that $|W(z)| \leq C \exp(C|z|^m)$ for all $z \in \mathbb{C}$.*

*Proof.* The first two conclusions follow directly from Lemma 5.5. The third requires a small additional argument in order to pass from upper bounds on $\psi_z^\pm$ at $\pm C|z|$, respectively, to bounds at 0. Assume $|z|$ to be large, Lemma 5.5 guarantees that $\psi_z^+$ and its first derivative are $O(\exp(C|z|^m))$ at $x = C|z|$, if $C$ is chosen to be sufficiently large. On the interval $[-C|z|, C|z|]$, the potential $Q^2(x, z)$ is $O(|z|^{2(m-1)})$. Therefore $|\partial_x^2 \psi_z^+| \leq C|z|^{2(m-1)} |\psi_z^+|$ on this interval. A simple comparison argument then gives $\psi_z^+ = O(\exp(C|z|^m))$ at $x = 0$, and likewise for its first derivative. The same reasoning applies to $\psi_z^-$, yielding the desired bound on the Wronskian. □

*Proof of Proposition 5.3.* By the first two conclusions of the preceding lemma, there exists an entire holomorphic function $W(z, \zeta)$ such that for each $\zeta \in \mathbb{C}^n$, $E_\zeta = \{z \in \mathbb{C} : W(z, \zeta) = 0\}$. Thus the exceptional set of all $\zeta$ for which there exist no nonlinear eigenvalues $z$ equals the set of all $\zeta$ for which $z \mapsto W(z, \zeta)$ has no zeros. By hypothesis, there exists at least one $\zeta_0$ for which $W(z, \zeta_0)$ has at least one zero. Therefore the set of all $\zeta \in \mathbb{C}^n$ for which there are no zeros is pluripolar [1]. □

## 6. Analysis of special cases

This section is devoted to the proof of Proposition 1.7. Let $Q(x, z) = x^{m-1} - z^k x^{m-k-1}$. $W(z)$ is an entire holomorphic function of $z^k$, so

$$\mathcal{W}(\zeta) = W(\zeta^{1/k})$$

is independent of the choice of $k$-th root, and hence is also entire holomorphic. By Corollary 5.7, $|\mathcal{W}(z)| \leq C \exp(C|z|^{m/k})$ for all $z \in \mathbb{C}$.

There are now two difficulties. Firstly, a superficial examination of [5],[8] might lead one to expect a lower bound $|W(z)| \geq c \exp(c|z|^m)$ for all $z \in \mathbb{R}$, but this is false[6] for $Q = x^2 + z^2$. Secondly, whereas the strategy in [5] and other closely related works was to show that $\mathcal{W}$ is an entire function of nonintegral order and hence must have zeros, it turns out that for the class of examples contemplated in Proposition 1.7, $\mathcal{W}$ is indeed an entire function of the expected order $m/k$, but this order may be integral. In that case a more refined analysis is required; this same point arose in the proof of Theorem 5.2 above.

Factor

$$Q(x, z) = x^{m-1-k}(x^k - z^k).$$

As in the proof of Theorem 5.2, the zeros of $\mathbb{R} \ni x \mapsto Q(x, z)$ are the key to the asymptotic behavior of $W(z)$, for $z^k \in \mathbb{R}$. Assume always that $z^k \in \mathbb{R}$. Then $x \mapsto Q(x, z)$ has two zeros when $k$ is odd, and has either one zero or three, depending on the sign of $z^k$, when $k$ is even.

The differential equation is $\partial_x^2 \psi_z^{\pm} = Q(x, z)^2 \psi_z^{\pm}$. Its solutions $\psi_z^{\pm}$ are real valued for real $x, z^k$. The asymptotics (5.3) imply that $\psi_z^+(x)$ is positive and $\partial_x \psi_z^+(x)$ negative for all sufficiently large $x \in \mathbb{R}^+$, and the differential equation then implies that $\psi_z^+$ is a positive, convex, decreasing function of $x$ on the whole real axis. Likewise $\psi_z^-$ is positive, convex, and increasing. Consequently

$$W(z) = \psi_z^+(0)\partial_x \psi_z^-(0) - \psi_z^-(0)\partial_x \psi_z^+(0) > \psi_z^+(0)\partial_x \psi_z^-(0). \qquad (6.1)$$

Thus in order to derive a lower bound for $W(z)$, it suffices to derive lower bounds for $\psi_z^+$ and for $\partial_x \psi_z^-$ at 0. We have $P(x, z) = m^{-1}x^m - (m-k)^{-1}z^k x^{m-k}$, so

$$P(z, z) = (m^{-1} - (m-k)^{-1})z^m = -c_0 z^m$$

where $c_0 > 0$.

---

[6]For $Q = x^2 + z^2$, $W$ is a polynomially bounded function on the real axis, as follows from the proof of Lemma 6.4 below, but satisfies the lower bound $c \exp(c|z|^3)$ on the imaginary axis (by the proof of Lemma 6.1).

**Lemma 6.1.** *There exists $c > 0$ such that*

$$\psi_z^+(z) \geq \exp(cz^m)$$

*for all sufficiently large $z \in \mathbb{R}^+$. If $k$ is even, so that $Q(-z, z) = 0$, then likewise*

$$\partial_x \psi_z^-(-z) \geq \exp(cz^m)$$

*for all sufficiently large $z \in \mathbb{R}^+$.*

Since $\psi_z^+(z) < \psi_z^+(0)$ and $\partial_x \psi_z^-(-z) < \partial_x \psi_z^-(0)$ by monotonicity, these same bounds hold at $x = 0$, as well. In particular, the lemma implies that $W(z) \geq \exp(cz^m)$ as $z \to +\infty$, when $k$ is even. More generally, it implies such a lower bound for $W$ whenever we can show that for some finite exponent $N$, $\partial_x \psi_z^-(0) \geq z^{-N}$ as $z \to +\infty$.

*Proof of Lemma 6.1.* The analysis of $\psi_z^-(-z)$ for even $k$ is completely parallel to the analysis for $\psi_z^+(z)$, so we treat only the latter. Set

$$\sigma = (m - 2)/2.$$

We will prove that for all sufficiently large $\rho \in \mathbb{R}^+$, for all $x \geq z + \rho z^{-\sigma}$, for $k = 0, 1$,

$$\partial_x^k \psi_z^+(x) = (-1)^k Q(x, z)^{k - \frac{1}{2}} e^{-P(x,z)}(1 + o(1)) \text{ as } \rho \to +\infty, \qquad (6.2)$$

uniformly in $x, z$ for all sufficiently large $z \in \mathbb{R}^+$. Since $|P(z + \rho z^{-\sigma}, z) - P(z, z)| \leq C_\rho$ uniformly for $z \geq 1$, since $Q(z + \rho z^{-\sigma}, z) \geq cz^{-N}$ for some finite $N$, and since $\psi_z^+$ is a decreasing function of $x$, this would imply that

$$\psi_z^+(z) \geq \psi_z^+(z + \rho z^{-\sigma}) \geq \exp(cz^m)$$

for some $c > 0$ as $z \to +\infty$, as desired.

To simplify notation in the proof, write $\psi = \psi_z^+$, $Q(x), P(x)$ for $Q(x, z), P(x, z)$, and write $f'$ for $\partial f / \partial x$. To prove (6.2) set

$$u = \begin{pmatrix} \psi \\ \psi' \end{pmatrix}, \qquad A = \begin{pmatrix} 0 & 1 \\ Q^2 & 0 \end{pmatrix} \qquad (6.3)$$

so that $u' = Au$. Set

$$S = \begin{pmatrix} 1 & 1 \\ Q - \frac{1}{2}Q'/Q & -Q - \frac{1}{2}Q'/Q \end{pmatrix} \qquad \text{for } x > z,$$

in which region $Q$ is always nonzero.

If $\rho$ is chosen to be sufficiently large, then for all sufficiently large $z$ and all $x \geq z + \rho z^{-\sigma}$, $0 < Q'(x, z)/Q^2(x, z) \leq \varepsilon(\rho)$ where $\varepsilon(\rho) \to 0$ as $\rho \to +\infty$; this calculation is contained in the proof of Lemma 6.2 below. We assume always that $x \geq z + \rho z^{-\sigma}$, where $\rho$ is sufficiently large. $S$ is then clearly invertible, and

$$S^{-1} = -\frac{1}{2}Q^{-1}\begin{pmatrix} -Q - \frac{1}{2}Q'/Q & -1 \\ -Q + \frac{1}{2}Q'/Q & 1 \end{pmatrix}.$$

This formal matrix calculation and others below may be found in [6], and the details will not be reproduced here.

The column vector $v$ satisfies the first order system of equations

$$v' = (B + E)v, \tag{6.4}$$

where

$$B = \begin{pmatrix} Q - \frac{1}{2}Q'/Q & 0 \\ 0 & -Q - \frac{1}{2}Q'/Q \end{pmatrix}$$

and $E$ is a continuous matrix valued function satisfying

$$|E| \leq C(Q')^2/|Q|^3 + C|Q''|/Q^2. \tag{6.5}$$

That is, each entry of $E$ satisfies this upper bound, uniformly in $x, z$ for $x \geq z + \rho z^{-\sigma}$ for any fixed $\rho$, assuming always that $z \geq 1$.

Set

$$w = \begin{pmatrix} 0 \\ Q^{-1/2}e^{-P} \end{pmatrix},$$

and

$$\Lambda(x) = \begin{pmatrix} Q^{-1/2}e^{P} & 0 \\ 0 & Q^{-1/2}e^{-P} \end{pmatrix}$$

so that $\Lambda' = B\Lambda$. Define the integral operator

$$Tf(x) = \int_x^\infty \Lambda(x)\Lambda^{-1}(y)(Ef)(y)\,dy, \tag{6.6}$$

acting on functions $f$ defined on $[z + \rho z^{-\sigma}, \infty)$ and taking values in the space of column vectors with two (real) entries.

If $v$ satisfies the integral equation

$$v(x) = w(x) - Tv(x) \tag{6.7}$$

then it satisfies the differential equation $v' = (B + E)v$. Indeed,

$$
\begin{aligned}
v'(x) = (w - Tv)'(x) &= w'(x) - \int_x^\infty \Lambda'(x)\Lambda^{-1}(y)(Ev)(y)\,dy + (Ev)(x) \\
&= (Bw)(x) - \int_x^\infty B(x)\Lambda(x)\Lambda^{-1}(y)(Ev)(y)\,dy + (Ev)(x) \\
&= Bw(x) - B(x)(Tv)(x) + (Ev)(x) \\
&= Bv(x) + Ev(x).
\end{aligned}
$$

Denote by $\varepsilon(\rho, z)$ any quantity that tends to zero as $\min(\rho, z) \to \infty$; this quantity is permitted to change from one line to the next.

**Lemma 6.2.** *The remainder coefficient matrix $E$ satisfies*

$$
\int_{z+z^{-\sigma}}^\infty |E(y)|\,dy \le C < \infty,
$$

*uniformly as $z \to +\infty$. Moreover*

$$
\int_{z+\rho z^{-\sigma}}^\infty |E(y)|\,dy \le \varepsilon(\rho, z).
$$

*Proof.* Write $A \sim B$ to mean that $A, B$ are positive quantities whose ratio is bounded above and below by positive constants independent of $z$, provided that $z$ is sufficiently large. Assume always that $x \ge z \gg 1$. Then

$$
\begin{aligned}
|Q(x, z)| &\sim x^{m-2}(x - z) \\
|Q'(x, z)| &\le C x^{m-2} \\
|Q''(x, z)| &\le C x^{m-3},
\end{aligned}
$$

where all derivatives indicated are taken with respect to $x$. Recalling the pointwise bound (6.5) for $E$, we obtain

$$
\int_{2z}^\infty |E(x, z)|\,dx \le C \int_{2z}^\infty \Big[ x^{2(m-2)-3(m-1)} + x^{m-3-2(m-1)} \Big]\,dx \le C z^{-m}.
$$

Likewise

$$
\begin{aligned}
\int_{z+\rho z^{-\sigma}}^{2z} |E(x, z)|\,dx &\le C \int_{z+\rho z^{-\sigma}}^{2z} \Big[ x^{2(m-2)-3(m-2)}(x-z)^{-3} + x^{m-3-2(m-2)}(x-z)^{-2} \Big]\,dx \\
&\le C z^{2-m} \int_{z+\rho z^{-\sigma}}^{2z} (x-z)^{-3}\,dx + C z^{1-m} \int_{z+\rho z^{-\sigma}}^{2z} (x-z)^{-2}\,dx \\
&\le C z^{2-m}(\rho z^{-\sigma})^{-2} + C z^{1-m}(\rho z^{-\sigma})^{-1} \\
&= C\rho^{-2} z^0 + C\rho^{-1} z^{1-m+\frac{m}{2}-1} \\
&\le \varepsilon(\rho, z).
\end{aligned}
$$

$\square$

For each $x \geq z + z^{-\sigma}$ introduce the norm

$$\|f\|_{*,x} = \sup_{y \geq x} Q^{1/2}(y) e^{P(y)} |f(y)| \qquad (6.8)$$

on the space of all continuous 2 dimensional column vector valued functions on $[x, \infty)$.

**Lemma 6.3.** *For each $\rho \geq 1$, for all $x \geq z + \rho z^{-\sigma}$,*

$$\|Tf\|_{*,x} \leq C_\rho \|f\|_{*,x}$$

*for all $f$ for which the right hand side is finite, uniformly for all $z \geq 1$ and uniformly in $x$. Moreover $C_\rho \to 0$ as $\rho \to +\infty$.*

*Proof.* Suppose that $\|f\|_{*,x} \leq 1$. Then

$$
\begin{aligned}
|Tf(x)| & \\
&\leq \int_x^\infty \left| \begin{pmatrix} Q^{-1/2}(x) e^{P(x)} Q^{1/2}(y) e^{-P(y)} & 0 \\ 0 & Q^{1/2}(x) e^{-P(x)} Q^{-1/2}(y) e^{P(y)} \end{pmatrix} \right| \\
&\qquad \cdot |E(y)| \, Q^{-1/2}(y) e^{-P(y)} \, dy \\
&\leq Q^{-1/2}(x) e^{-P(x)} \int_x^\infty |E(y)| dy \\
&\leq C Q^{-1/2}(x) e^{-P(x)}
\end{aligned}
$$

by Lemma 6.2, because $e^{-P}$ and $Q^{-1/2}$ are both decreasing functions of $x$ for $x \geq z$. $\qquad \square$

If $\rho$ is chosen to be sufficiently large, then we find that the map $f \mapsto w - Tf$ is a strict contraction on the space of all continuous functions on $[z + \rho z^{-\sigma}, \infty)$ for which $\|f\|_{*,[z+\rho z^{-\sigma}]}$ is finite. Therefore by the contraction mapping principle, there exists a unique solution of the integral equation (6.7) for which this norm is finite. We define $v$ to be that solution.

It follows from the same reasoning together with the second conclusion of Lemma 6.2 that

$$|(v-w)(y)| \leq \varepsilon(\rho) Q^{-1/2}(y) e^{-P(y)} \qquad (6.9)$$

for all $y \geq z + \rho z^{-\sigma}$, where $\varepsilon(\rho) \to 0$ as $\rho \to \infty$.

Defining $u = Sv$, $u$ satisfies the differential equation $u' = Au$. Thus

$$u = Sv = \begin{pmatrix} Q^{-1/2} e^{-P} (1 + \varepsilon(\rho)) \\ Q^{+1/2} e^{-P} (1 + \varepsilon(\rho)) \end{pmatrix}$$

for all $x \geq z + \rho z^{-\sigma}$, since $|Q'|/Q^2 \to 0$ as $\rho \to \infty$ in that region, uniformly for $z \geq 1$.

$\psi_z^+$ is defined (for large $z \in \mathbb{R}^+$) in terms of $u$ by equation (6.3). In particular,

$$\psi_z^+(x) = Q^{-1/2}(x, z)e^{-P(x,z)}(1 + \varepsilon(\rho))$$

for $x \geq z + \rho z^{-\sigma}$. At $x = z + \rho z^{-\sigma}$, this is bounded below by $e^{cz^m}$ for some $c > 0$, as $z \to +\infty$. This concludes the proof of Lemma 6.1. $\qquad \square$

**Lemma 6.4.** *If $k$ is odd then there exist $C, N \in \mathbb{R}^+$ such that for all $z \geq 1$,*

$$C^{-1}z^{-N} \leq \psi_z^-(0) \leq Cz^N,$$

*and the same holds for $\partial_x \psi_z^-(0)$.*

*Proof.* Consider first the case where $m$ is even. Let $B, E, W, S, \Lambda$ be as defined in the proof of Lemma 6.1. Then $v' = (B + E)v$. Redefine

$$Tf(x) = \int_{-\infty}^x \Lambda(x)\Lambda^{-1}(y)(Ev)(y)\,dy$$

and denote by $v$ a putative solution of the integral equation

$$v = w + Tv.$$

Redefine

$$\sigma = k/(m - k).$$

We again have

$$\int_{-\infty}^{-\rho z^{-\sigma}} |E(x, z)|\,dx \leq C\varepsilon(\rho, z).$$

Indeed, for $z \geq 1$ and $x \leq 0$,

$$|Q(x, z)| \sim |x|^{m-k-1}(|x|^k + z^k)$$
$$|Q'(x, z)| \leq C|x|^{m-k-2}(|x|^k + z^k)$$
$$|Q''(x, z)| \leq C|x|^{m-k-3}(|x|^k + z^k).$$

Thus

$$\int_{-\infty}^{-z} |E(x, z)|\,dx \leq C \int_{-\infty}^{-z} \Big[|x|^{2(m-2)-3(m-1)} + |x|^{m-3-2(m-1)}\Big]\,dx \leq Cz^{-m},$$

while the integral of $|E(x, z)|$ over $[-z, -\rho z^{-\sigma}]$ is majorized by

$$C\int_{-z}^{-\rho z^{-\sigma}} \Big[z^{2k-3k}|x|^{2(m-k-2)-3(m-k-1)} + z^{k-2k}|x|^{m-k-3-2(m-k-1)}\Big]\,dx$$

$$\leq Cz^{-k}\int_{-z}^{-\rho z^{-\sigma}} |x|^{-m+k-1}\,dx$$

$$\leq Cz^{-k}(\rho z^{-\sigma})^{k-m}$$

$$= C\rho^{k-m}.$$

Redefine

$$\|f\|_{*,x} = \sup_{y \le x} |Q|^{1/2}(y) e^{P(y)} |f(y)|.$$

Because $k$ is odd, $|Q|^{-1/2} e^{-P}$ is a monotone increasing function on $(-\infty, 0]$. Therefore for $x \le 0$, for any $f$ with finite $*, x$ norm,

$$|Tf(x)| \le \|f\|_{*,x} \int_{-\infty}^{x} \left[ |Q|^{-1/2}(y) e^{-P(y)} + |Q|^{1/2}(x) e^{P(x)} |Q|^{-1}(y) e^{-2P(y)} \right] |E(y)| \, dy$$

$$\le \|f\|_{*,x} |Q|^{-1/2}(x) e^{-P(x)} \int_{-\infty}^{x} |E(y)| \, dy.$$

Thus for $x \le -\rho z^{-\sigma}$,

$$\|Tf\|_{*,x} \le \varepsilon(\rho, z) \|f\|_{*,x}.$$

Consequently, when $\rho, z$ are both sufficiently large, $T$ is a strict contraction on the space of all continuous two dimensional column vector valued functions $f$ having finite $*, -\rho z^{-\sigma}$ norms, and so there exists a unique solution with finite norm, $v$, of the integral equation $v = w + Tv$. Moreover $\|v - w\|_{*,-\rho z^{-\sigma}} \le \varepsilon(\rho, z)$.

Defining $u = Sv$, we have

$$u = \begin{pmatrix} \psi_z^- \\ \partial_x \psi_z^- \end{pmatrix}.$$

The equation $v = w + Tv$ and bound on $Tv = v - w$ bound obtained in the preceding paragraph imply that

$$u = \begin{pmatrix} |Q|^{-1/2} e^{-P} \cdot (1 + \varepsilon(\rho, z)) \\ |Q|^{+1/2} e^{-P} \cdot (1 + \varepsilon(\rho, z)) \end{pmatrix}$$

for $x \le -\rho z^{-\sigma}$. At $x = -\rho z^{-\sigma}$, $|Q| \sim z^k |x|^{m-k-1} \sim z^{k-[(m-k-1)k/(m-k)]} = z^{k/(m-k)} = z^{\sigma}$. For $-\rho z^{-\sigma} \le x \le 0$ we have $Q(x) = O(z^{\sigma})$, so for all such $x$, $|P(x, z)| = |\int_0^x Q(y, z) dy| = O(1)$. Thus

$$C^{-1} z^{-\sigma/2} \le \quad \psi_z^-(-\rho z^{-\sigma}) \le C z^{-\sigma/2},$$

$$C^{-1} z^{+\sigma/2} \le \partial_x \psi_z^-(-\rho z^{-\sigma}) \le C z^{+\sigma/2}.$$

For $x \in [-\rho z^{-\sigma}, 0]$ we have a differential inequality

$$|(\psi_z^-)''(x)| \le C z^{2\sigma} \psi_z^-(x).$$

From a simple comparison argument and the upper bounds for $\psi_z^-$ and its first derivative at $-\rho z^{-\sigma}$, it follows that $\psi_z^-(0) = O(z^{-\sigma/2})$ and $\partial_x \psi_z^-(0) = O(z^{\sigma/2})$. On the other hand, convexity of $\psi_z^-$ implies that

$$\partial_x \psi_z^-(0) > \partial_x \psi_z^-(-\rho z^{-\sigma}) \ge C^{-1} z^{\sigma/2}.$$

The case of odd $m$ is treated in the same way, except that the formulas are changed as in Lemma 5.5, with $P$ replaced by $-P$ throughout. □

**Corollary 6.5.** *For every $m \geq 2$ and $k \geq 1$, $W$ is an entire holomorphic function of order exactly $m$.*

*Proof.* If $k$ is even then

$$\psi_z^+(0) \geq \psi_z^+(z) \geq e^{cz^m}, \tag{6.10}$$

while

$$\partial_x \psi_z^-(0) \geq \partial_x \psi_z^-(-z) \geq e^{cz^m}$$

and hence

$$W(z) > e^{cz^m}$$

for all large $z \in \mathbb{R}^+$. Since $|W(z)| \leq C \exp(C|z|^m)$ for all $z \in \mathbb{C}$, $W$ has order exactly $m$.

If $k$ is odd then (6.10) remains valid, while the lower bound of Lemma 6.4 implies that as $\mathbb{R} \ni z \to +\infty$,

$$W(z) > \psi_z^+(0)\partial_x \psi_z^-(0) \geq C^{-1} z^{-N} e^{cz^m} > e^{\delta z^m}$$

for some $\delta > 0$. □

**Lemma 6.6.** *If $k$ is even then there exist $C, N \in \mathbb{R}^+$ such that whenever $z^k \in \mathbb{R}^-$ and $|z|$ sufficiently large,*

$$C^{-1}|z|^{-N} \leq |W(z)| \leq C|z|^N.$$

*Proof.* For $x \in \mathbb{R}$ and $z^k \in \mathbb{R}^-$, for $k$ even, $|Q(x, z)| \sim |x|^{m-k-1}(|x|^k + |z|^k)$. Both $\psi_z^\pm$ may therefore be treated as was $\psi_z^-$ in Lemma 6.4, for $x \geq 0$ in the case of the plus sign, and for $x \leq 0$ in the case of the minus sign. □

**Lemma 6.7.** *For any $m \geq 2$ and $k \geq 1$,*

$$|W(z)| \leq C|z|^N$$

*for all $|z| \geq 1$ such that $z^k \in i\mathbb{R}$.*

*Proof.* When $z^k \in i\mathbb{R}$ and $x \in \mathbb{R}$, $|Q(x, z)| \sim |x|^{m-k-1}(|x|^k + |z|^k)$. Now $|e^{-P(x)}| \equiv e^{-x^m/m}$, so the $*, x$ norms are slightly different though they are defined by the same formal expressions. $|e^{-P}|$ and $|Q|$ are still monotone. The proof is thus nearly identical to that of Lemma 6.4. □

With these lemmas in hand we now prove Proposition 1.7. Any entire holomorphic function $F$ of finite, strictly positive order either has infinitely many zeros, or takes the form $\exp(R)$ for some polynomial $R$ whose degree equals the order of $F$. From Corollary 6.5 we know that in all cases of Proposition 1.7, $\mathcal{W}$ is an entire function of order exactly $m/k$. Hence if $m/k$ is not an integer, $\mathcal{W}$ must have zeros.

Note that whenever $k$ is odd, $\mathcal{W}$ is an even function, since $\partial_x^2$ and $(x^{m-1} - \zeta x^{m-k-1})^2$ are invariant under the substitution $(x, \zeta) \mapsto (-x, -\zeta)$. Consequently if $\mathcal{W}$ had no zeros and hence equalled $\exp(R)$ for some polynomial, $R$ would also have to be even, and hence $R$ would have to have even degree, so $m/k$ would be forced to be even.

Suppose $m$ is odd. If $k$ is even, then $m/k$ is not integral, so $\mathcal{W}$ has zeros. If $k$ is odd, then $m/k$ is odd, and hence by the preceding paragraph, $\mathcal{W}$ must have zeros.

Suppose $k$ is even. By Lemma 6.6, $|\mathcal{W}(\zeta)|$ and its reciprocal have at most polynomial growth as $\mathbb{R}^- \ni \zeta \to -\infty$, while $|\mathcal{W}(\zeta)| \geq \exp(c\zeta^{m/k})$ as $\mathbb{R}^+ \ni \zeta \to +\infty$. There exists no polynomial $R$ of degree $m/k$ such that $\exp(R(\zeta))$ behaves in this way, so $\mathcal{W}$ cannot have a polynomial logarithm, hence must have zeros.

Suppose $m$ is divisible by 4. By the preceding paragraph it suffices to examine the case where $k$ is odd. Then the order $m/k$ of $\mathcal{W}$ is also divisible by 4. Suppose $\mathcal{W} = \exp(R)$ where $R$ is a polynomial. Then since $|\mathcal{W}(\zeta)| \geq \exp(c|\zeta|^{m/k})$ for large $\zeta \in \mathbb{R}^+$, we have $\operatorname{Re}(R(\zeta)) \geq c|\zeta|^{m/k}$ for large $\zeta \in \mathbb{R}^+$, for some $c > 0$. Since the degree $m/k$ of $R$ is divisible by 4, this forces $\operatorname{Re}(R(\zeta)) \geq |\zeta|^{m/k}/2$ whenever $\zeta \in i\mathbb{R}$ has sufficiently large modulus. This contradicts Lemma 6.7. The proof of Proposition 1.7 is therefore complete. $\qquad\square$

## 7. Sufficiency

In this section we outline a proof of Theorem 1.2. The method is a straightforward combination of the FBI transform with a subset of the machinery used in Section 4 to prove the negative results.

Define an FBI transform of any compactly supported function by

$$\mathcal{F}u(x, t; \xi, \tau) = \int_{\mathbb{R}^2} e^{i[(x-y)\xi + (t-s)\tau] - |\tau|(x-y)^2 - |\tau|(t-s)^2} u(y, s) \, dy \, ds. \qquad (7.1)$$

It will always be assumed that $|\xi| \leq |\tau|$. For any point $p = (x_0, t_0; \xi_0, \tau_0) \in T^*\mathbb{R}^2$ with $|\xi_0| < |\tau_0|$, $p$ belongs to the complement of $WF_a(u)$, the analytic wave front set of $u$, if and only if there exist $\varepsilon > 0$ and a conic neighborhood $\Gamma$ of $p$ such that

$$|\mathcal{F}u(x, t; \xi, \tau)| \leq C \exp(-\varepsilon|\tau|) \qquad (7.2)$$

for all $(x, t; \xi, \tau) \in \Gamma$ [17]. Fix a coordinate system in which the span of $\{X_1, X_2\}$ coincides with the span of $\{\partial_x, x^{m-1}\partial_t\}$. Suppose that $Lu \in C^\omega$ in some neighborhood $U$ of 0; we wish to prove that $u \in C^\omega$ in the intersection of $U$ with a fixed neighborhood $U_0$ of 0. Analyticity at other points follows by ellipticity where $x \neq 0$, and by a change of coordinates to reduce the case $x = 0$ to $x = t = 0$. Since $L$ is $C^\infty$ hypoelliptic, $u \in C^\infty$ near 0.

The analytic wave front set of $u$ is contained in the union of the characteristic variety $\Sigma$ of $L$ with the analytic wave front set of $Lu$ [17],[13]. In our coordinates, $\Sigma$ is simply $\{(x, t, \xi, \tau) : x = \xi = 0\}$. Thus in order to prove that $u \in C^\omega(U \cap U_0)$, it suffices to prove that (7.2) holds for all $(x', t')$ in any fixed compact subset of $U \cap U_0$ and for all $(\xi, \tau)$ satisfying $|\xi| \leq |\tau|$.

Fix any $(\xi, \tau)$ satisfying $|\xi| \leq |\tau|$ and a point $(x', t')$ near 0. Define a differential operator $L_\tau$, acting with respect to the variables $(x, t)$ and depending on the parameters $(t', \tau)$, by

$$L_\tau = (e_{(t', \tau)})^{-1} \circ L^* \circ e_{(t', \tau)}$$

where $e_{(t', \tau)}(x, t) = \exp(i(t' - t)\tau - |\tau|(t' - t)^2)$ and $L^*$ denotes the transpose of $L$. The notation $L_\tau$ is misleading because the operator depends on $t'$ as well, but the dependence on $\tau$ is of greater importance.

The following lemma is analogous to Lemma 4.1. Denote by $D_\delta$ the open disk in $\mathbb{C}$ of radius $\delta$, centered at the origin.

**Lemma 7.1.** *For any sufficiently small $\delta > 0$ there exist $C, \varepsilon \in \mathbb{R}^+$ such that for each $(x', t') \in \mathbb{R}^2$ satisfying $|(x', t')| \leq \delta/8$ and for each $(\xi, \tau) \in \mathbb{R}^2$ satisfying $|\xi| \leq |\tau|$, there exists a function $g = g_{(x', t', \xi, \tau)}(x, t)$ defined for $(x, t) \in (-\delta, \delta) \times D_\delta$ and holomorphic with respect to $t$, satisfying the equation*

$$L_\tau g(x, t) = e^{i(x' - x)\xi - |\tau|(x' - x)^2} + O(e^{-\varepsilon|\tau|}) \tag{7.3}$$

*with the bounds*

$$|g(x, t)| \leq C e^{-\varepsilon|\tau|} \qquad \text{for all } |x| \geq \delta/4 \text{ and } t \in D_\delta \tag{7.4}$$

$$|g(x, t)| \leq C \qquad \text{for all } (x, t) \in [-\delta, \delta] \times D_\delta. \tag{7.5}$$

*Proof of Theorem 1.2.* Granting the lemma, and supposing that $Lu \in C^\omega$ near 0, there exist $\delta$ and a $C^\infty$ function $v$ supported in $(-\delta, \delta) \times (-2\delta, 2\delta)$, such that $Lv \equiv Lu$ for all $(x, t) \in \mathbb{R}^2$ satisfying $|x| \leq \delta/2$ and $|t| \leq \delta$. Consider any $(\xi, \tau)$ satisfying $|\xi| \leq |\tau|$, and suppose for simplicity of notation that $\tau > 0$. Consider likewise any

$(x', t') \in \mathbb{R}^2$ sufficiently close to the origin. Then writing $g = g_{(x', t', \xi, \tau)}$,

$$
\begin{aligned}
\mathcal{F}v(x', t'; \xi, \tau) &= \iint v(x, t) e^{i(x'-x)\xi - \tau(x'-x)^2} e^{i(t'-t)\tau - \tau(t'-t)^2} \, dx \, dt \\
&= \iint v e^{i(t'-t)\tau - \tau(t'-t)^2} L_\tau g \quad + O(e^{-\varepsilon\tau}) \\
&= \iint v \cdot L^* \big( e^{i(t'-t)\tau - \tau(t'-t)^2} g \big) \quad + O(e^{-\varepsilon\tau}) \\
&= \iint Lv \cdot e^{i(t'-t)\tau - \tau(t'-t)^2} g \quad + O(e^{-\varepsilon\tau}).
\end{aligned}
$$

Replace the two dimensional contour of integration by the contour in $\mathbb{R} \times \mathbb{C}$ parametrized by $(x, t) \mapsto (x, t + ih(x, t))$ where $h$ is smooth, takes values in $(-1/2, 0]$, is identically equal to zero where $|(x, t)| \geq \delta/2$, and is strictly negative but small where $|(x, t)| < \delta/2$. Since $Lv$ is a holomorphic function of $t$ where $|(x, t)| < \delta/2$, we obtain

$$
\begin{aligned}
\mathcal{F}v(x', t'; \xi, \tau) &= O(e^{-\varepsilon\tau}) \\
&+ \iint Lv(x, t + ih(x, t)) e^{i(t'-t-ih(x,t))\tau - \tau(t'-t-ih(x,t))^2} g(x, t + ih(x, t)) \, J(x, t) \, dx \, dt
\end{aligned}
$$

where $J(x, t) = 1 + i\partial h/\partial t$.

Where $|(x, t)| \leq \delta/2$, both $Lv(x, t_i h)$ and $g(x, t + ih)$ are bounded uniformly in $(x', t'; \xi, \tau)$, while the real part of $i(t' - t - ih(x, t))\tau - \tau(t' - t - ih(x, t))^2$ equals $[h(x, t) - h(x, t)^2]\tau$, which is $\leq -\varepsilon\tau$ for some $\varepsilon > 0$, provided that $(x', t')$ lies in a sufficiently small neighborhood of $0$ that is independent of $\xi, \tau$. Thus the exponential factor in the integral is $O(\exp(-\varepsilon\tau))$, hence so is the integrand.

Where $|x| \geq \delta/4$, both $Lv$ and the exponential factor are still $O(1)$, and $g$ is $O(\exp(-\varepsilon\tau))$. Lastly, where $|t| \geq \delta/4$, the exponential factor is $O(\exp(-\varepsilon\tau))$ and $Lv, g$ are $O(1)$. Thus $\mathcal{F}v(x', t'; \xi, \tau) = O(\exp(-\varepsilon\tau))$, as desired.  $\square$

The proof of Lemma 7.1 is parallel to that of Lemma 4.1, with one essential difference: the real part of $\tilde{\Theta}(x, t)^2$ is nonnegative for all $(x, t)$ in some neighborhood of $0$ in $\mathbb{R} \times \mathbb{C}$ if (and only if) $\Theta(x, t) \equiv x^{m-1}$. Consequently the proof of Lemma 4.1 may be executed in a fixed neighborhood of $0$ with respect to the $t$ coordinate, rather than in a neighborhood that shrinks to $\{0\}$ as $\tau \to \infty$. This in turn allows us to carry out the Neumann series to $c\tau$ terms, rather than to $c\tau^{p/m}$ terms, resulting in an error that is $O(\exp(-\varepsilon\tau))$, rather than merely $O(\exp(-\varepsilon\tau^{p/m}))$ where $p/m = 1 - (mq)^{-1}$.

Denote by $D$ any open disk in $\mathbb{C}$ centered at the origin. The Sobolev spaces appropriate for the Neumann series argument in the present context are as follows:

$$\|f\|^2_{\mathcal{H}^0_{\rho,\tau}(\mathbb{R}\times D)} = \iint_{\mathbb{R}\times D} |f(x)|^2 \, (1+\tau^{2/m}x^2)^{-(m-1)} \, e^{\tau\rho\eta(x)} \, dx \, dt d\bar{t}$$

$$\|f\|^2_{\mathcal{H}^1_{\rho,\tau}(\mathbb{R}\times D)} = \iint_{\mathbb{R}\times D} \left( (1+\tau^{2/m}x^2)^{-(m-1)}|\partial_x f|^2 + \tau^{2/m}|f|^2 \right) e^{\tau\rho\eta(x)} \, dx \, dt d\bar{t}$$

$$\|f\|^2_{\mathcal{H}^2_{\rho,\tau}(\mathbb{R}\times D)} = \iint_{\mathbb{R}\times D} \Big( (1+\tau^{2/m}x^2)^{-(m-1)}|\partial_x^2 f|^2 + \tau^{2/m}|\partial_x f|^2$$
$$+ \tau^{4/m}(1+\tau^{2/m}x^2)^{(m-1)}|f|^2 \Big) e^{\tau\rho\eta(x)} \, dx \, dt d\bar{t},$$

where $\eta \in C_0^\infty$ is nonnegative, $\eta(x) \equiv 0$ for all $|x| \leq \delta/8$ and $\eta(x) > 0$ for all $|x| \geq \delta/4$. The parameter $\rho$ is to be chosen to be positive but sufficiently small, and is independent of $\tau$.

In order to work in these spaces we must modify $L$ so that its coefficients are defined for all $(x,t)$ in the Cartesian product of $\mathbb{R}$ with a complex neighborhood of 0, and are holomorphic with respect to $t$. $\tilde{\Theta}(x,t)$ is identically equal to $x^{m-1}$ near 0, so we extend it to equal $x^{m-1}$ everywhere. $L$ takes the form $(a_{1,1}\partial_x + a_{1,2}\tilde{\Theta}\partial_t)^2 + (a_{2,1}\partial_x + a_{2,2}\tilde{\Theta}\partial_t)^2$ near 0, so it suffices to extend the coefficients $a_{i,j}$ to $\mathbb{R}\times\mathbb{C}$ so that they are independent of $x$ outside a small neighborhood of the origin, and so that the coefficient matrix $(a_{i,j})$ is real and invertible everywhere on the product of $\mathbb{R}$ with a small real neighborhood of 0, and depends holomorphically on $t$ in $\mathbb{R} \times D_\delta$ for some $\delta > 0$.

After conjugation, $\partial_t$ takes the form

$$e_{(t',\tau)}^{-1} \circ \partial_t \circ e_{(t',\tau)} = \partial_t - i\tau + 2(t'-t)\tau = -i\tau(1 + 2i(t'-t)) + \partial_t.$$

Write

$$L^* = P(x,t,\partial_x,\tilde{\Theta}\partial_t) = (a_{1,1}\partial_x + a_{1,2}\tilde{\Theta}\partial_t)^2 + (a_{2,1}\partial_x + a_{2,2}\tilde{\Theta}\partial_t)^2 + b_1\partial_x + b_2\tilde{\Theta}\partial_t + b_3$$

where the coefficients $a_{i,j}$ and $b_i$ are real and analytic. Then the conjugated operator $L_\tau$ takes the form

$$L_\tau = P(x,t,\partial_x,[-i\tau(1 + 2i(t'-t)) + \partial_t]\tilde{\Theta}).$$

Define the ordinary differential operators

$$A_t = P(x,t,\partial_x,[-i\tau(1 + 2i(t'-t))]\tilde{\Theta}).$$

which depend on the parameters $t,t',\tau$. Define an operator $\mathcal{A}$, acting on functions of $(x,t) \in \mathbb{R} \times D$ that are holomorphic with respect to $t$ in a small disk $D \subset \mathbb{C}$ centered at 0, by letting $A_t$ act with respect to the variable $x$, for each $t$.

**Lemma 7.2.** *For all $t, t' \in \mathbb{C}$ and $\rho \in \mathbb{R}$ sufficiently close to $0$ and all sufficiently large $\tau \in \mathbb{R}^+$, the map*

$$A_t : \mathcal{H}^2_{\rho,\tau}(\mathbb{R}) \mapsto \mathcal{H}^0_{\rho,\tau}(\mathbb{R})$$

*is invertible, uniformly in $t, t', \tau$. For all sufficiently small $t', \delta, \rho$ and sufficiently large $\tau$, the map*

$$\mathcal{A} : \mathcal{H}^2_{\rho,\tau}(\mathbb{R} \times D_\delta) \mapsto \mathcal{H}^0_{\rho,\tau}(\mathbb{R} \times D_\delta)$$

*is invertible, uniformly in $t', \delta, \tau$.*

*Proof.* The main point is the inequality

$$\|\phi\|^2_{\mathcal{H}^2_{0,\tau}(\mathbb{R})} \le C\|A_t\phi\|^2_{\mathcal{H}^0_{0,\tau}(\mathbb{R})} \qquad \text{for all } \phi \in C^2_0(\mathbb{R}). \tag{7.6}$$

To prove this, consider

$$-\operatorname{Re} \int_{\mathbb{R}} A_t\phi \cdot \bar{\phi} \ge c \int_{\mathbb{R}} |\partial_x\phi|^2 + \tau^2 x^{2(m-1)}|\phi|^2 - C \int_{\mathbb{R}} |\phi| \cdot \left( |\partial_x\phi| + \tau|x^{m-1}\phi| + |\phi| \right)$$

for some $c, C \in \mathbb{R}^+$, provided that $|t - t'|$ is sufficiently small. Since

$$\tau^{2/m} \int_{\mathbb{R}} |\phi|^2 \le C \int_{\mathbb{R}} |\partial_x\phi|^2 + \tau^2 x^{2(m-1)}|\phi|^2,$$

the last inequality plus Cauchy-Schwarz lead to

$$-\operatorname{Re} \int_{\mathbb{R}} A_t\phi \cdot \bar{\phi} \ge c' \int_{\mathbb{R}} |\partial_x\phi|^2 + \tau^{2/m}\left(1 + \tau^{1/m}|x|\right)^{2(m-1)}|\phi|^2$$

for some $c' > 0$, provided that $\tau$ is sufficiently large. Applying Cauchy-Schwarz to the left hand side yields

$$\int_{\mathbb{R}} |\partial_x\phi|^2 + \tau^{2/m}\left(1 + \tau^{1/m}|x|\right)^{2(m-1)}|\phi|^2 \le C \int_{\mathbb{R}} |A_t\phi|^2 \tau^{-2/m}\left(1 + \tau^{1/m}|x|\right)^{-2(m-1)}$$
$$\le C\tau^{2/m}\|A_t\phi\|^2_{\mathcal{H}^0_{0,\tau}}. \tag{7.7}$$

In order to prove (7.6), it remains only to bound the $\mathcal{H}^0_{0,\tau}$ norm of $\partial_x^2\phi$ by the $\mathcal{H}^0_{0,\tau}$ norm of $A_t\phi$. The operator $A_t$ is a quadratic polynomial, with $C^\omega$ coefficients, in $\partial_x$ and $\tau\tilde{\Theta}$. The result of applying any monomial of degree two or less, excepting only $\partial_x^2$, to $\phi$ has already been estimated in the $\mathcal{H}^0_{0,\tau}$ norm in (7.7). The missing term $\partial_x^2$ may be expressed as a linear combination, with $C^\omega$ coefficients, of $A_t$ and all the other monomials. Therefore the $\mathcal{H}^0_{0,\tau}$ norm of $\partial_x^2\phi$ is majorized by

$$\|\partial_x^2\phi\|^2_{H^0_{0,\tau}} \sim \int_{\mathbb{R}} (1 + \tau^{1/m}|x|)^{-2(m-1)}|\partial_x^2\phi|^2 \, dx$$
$$\le C\|A_t\phi\|^2_{H^0_{0,\tau}} + C \int_{\mathbb{R}} \tau^{2/m}|\partial_x\phi|^2 + C\tau^{4/m}\left(1 + \tau^{1/m}|x|\right)^{2(m-1)}|\phi|^2$$
$$\le C\|A_t\phi\|^2_{H^0_{0,\tau}}.$$

Thus we obtain the inequality (7.6) for $\rho = 0$. The case of small $|\rho|$ follows from the case $\rho = 0$ by conjugating with $e^{\rho \tau \eta(x)}$, as in the final paragraph of the proof of Lemma 3.1 of [3].

The spaces $\mathcal{H}^k_{\rho,\tau}$ are defined so that $A_t$ is automatically a bounded operator from $\mathcal{H}^2$ to $\mathcal{H}^0$. Inequality (7.6) implies that $A_t$ has closed range; the same analysis applies to its transpose and consequently invertibility of $A_t$ follows from a duality argument as in the proof of Lemma 3.1 of [3]. Since $A_t$ depends holomorphically on $t$ and is invertible uniformly in $t$, $A_t^{-1}$ also depends holomorphically on $t$. Therefore as an operator from $\mathcal{H}^2_{\rho,\tau}(\mathbb{R} \times D_\delta)$ to $\mathcal{H}^0_{\rho,\tau}(\mathbb{R} \times D_\delta)$, $\mathcal{A}$ may be inverted by applying the inverse of $A_t$ in the $x$ variable for each value of $t$. $\qquad\square$

With Lemma 7.2 in hand, imitating the proof of Lemma 4.1 leads directly to the conclusions of Lemma 7.1, hence to Theorem 1.2. $\qquad\square$

## 8. Microlocal analytic hypoellipticity for factored operators

Let $X_1, X_2$ be vector fields defined in a neighborhood of $p = 0$, satisfying the hypotheses of Theorem 1.8. Assume that $X_1, X_2$ are linearly dependent at the point 0. Suppose that $L = (X_1 + iX_2)(X_1 - iX_2) + c_1 X_1 + c_2 X_2 + c_3$ where $c_j \in C^\omega$ are complex valued.

Let $(x, t)$ be a system of coordinates with the properties of Lemma 1.3. Writing $X_1 = a_{1,1} \partial_x + a_{1,2} \tilde{\Theta} \partial_t$ and $X_2 = a_{2,1} \partial_x + a_{2,2} \tilde{\Theta} \partial_t$,

$$[X_1, X_2] = \text{Det}(a)(\partial_x \tilde{\Theta}) \partial_t + O(X_1, X_2)$$

where $\text{Det}(a)$ is the determinant of the matrix $a = (a_{i,j})$. This determinant is real valued and nonvanishing near 0 in $\mathbb{R}^2$; by replacing $t$ by $-t$ if necessary, we may require it to be negative.

Let $(\xi, \tau)$ be coordinates dual to $(x, t)$. Near 0, decompose the characteristic variety $\Sigma$ of $L$ as the disjoint union $\Sigma^+ \cup \Sigma^-$, where $\Sigma^+ = \Sigma \cap \{\tau > 0\}$. Note that the fiber of $\Sigma$ over 0 is $\{(\xi, \tau) : \xi = 0\}$.

**Lemma 8.1.** *Assume that $\{X_1, X_2\}$ is a pseudoconvex pair near 0. Then the function $\partial_x \Theta$ is everywhere nonnegative, in some neighborhood of 0, and the type $m$ at 0 is even.*

*Proof.* The pseudoconvexity hypothesis is

$$\partial_x \tilde{\Theta} = h + b\tilde{\Theta}, \tag{8.1}$$

where $h, b$ are smooth functions and $h$ does not change sign. Consider first the case where $h \geq 0$.

Expand $\partial_x \tilde{\Theta}$ and $\tilde{\Theta}$ in Taylor series, and assign weights $1, q^{-1}$ respectively to $x, t$, as in Section 2. Only monomials with weights $\geq m-2$ arise in the expansion of $\partial_x \tilde{\Theta}$, and the sum of all monomials having weight $m-2$ is equal to $\partial_x \Theta$. Any smooth multiple of $\tilde{\Theta}$ itself involves only terms of weights $m-1$ and greater. Thus $h = \partial_x \Theta$ modulo terms of weight $> (m-2)$. Consequently for any $(y, s)$,

$$\begin{aligned}
\partial_x \Theta(y, s) &= \lim_{\varepsilon \to 0} \varepsilon^{-(m-2)} \partial_x \tilde{\Theta}(\varepsilon y, \varepsilon^{1/q} s) \\
&= \lim_{\varepsilon \to 0} \varepsilon^{-(m-2)} \partial_x \tilde{\Theta}(\varepsilon y, \varepsilon^{1/q} s) - \varepsilon^{-(m-2)} b(\varepsilon y, \varepsilon^{1/q} s) \tilde{\theta}(\varepsilon y, \varepsilon^{1/q} s) \\
&= \lim_{\varepsilon \to 0} \varepsilon^{-(m-2)} h(\varepsilon y, \varepsilon^{1/q} s)
\end{aligned}$$

is a limit of nonnegative quantities, hence must be nonnegative.

In the case $h \leq 0$, the same reasoning leads to the conclusion that $\partial_x \Theta(x, t) \leq 0$ for all $(x, t) \in \mathbb{R}^2$. This is impossible, since $\partial_x \Theta(x, 0) \equiv (m-1) x^{m-2}$. Likewise nonnegativity of $\partial_x \Theta$ forces $m$ to be even, for otherwise $(m-1) x^{m-2}$ would change sign. $\qquad \square$

The principal symbol of the commutator satisfies

$$\sigma_1(i[X_1, X_2]) = -\operatorname{Det}(a)(\partial_x \tilde{\Theta})\tau + O(\sigma_1(iX_1), \sigma_1(iX_2)).$$

Thus modulo the span of the symbols of $X_1, X_2$, this is nonnegative when $\tau \geq 0$. The function $\tilde{\Theta}$ cannot vanish identically along any line segment where the $t$ coordinate is constant, by the bracket hypothesis. Therefore arbitrarily close to any point where $\tilde{\Theta}(x, t) = 0$ there exist points where $\tilde{\Theta} \neq 0$ but the ratio $|\partial_x \tilde{\Theta}/\tilde{\Theta}|$ is arbitrarily large. For $(x, t, \xi, \tau) \in \Sigma^-$ one has $\tau < 0$, so by the pseudoconvexity hypothesis, $-\tau \partial_x \tilde{\Theta}(x, t) \geq 0$ modulo a bounded multiple of $|\tau \tilde{\Theta}(x, t)|$, so $\tau \partial_x \tilde{\Theta}(x, t)$ must be negative, even modulo any bounded multiple of $|\tau \tilde{\Theta}(x, t)|$, wherever the ratio $|\partial_x \tilde{\Theta}/\tilde{\Theta}|$ is sufficiently large. Consequently $\Sigma^-$ has no conic neighborhood in which $\tau \partial_x \tilde{\Theta}$ is nonpositive modulo the span of the symbols of $X_1, X_2$, so the half line bundle $\Sigma^+$ is uniquely specified by the existence of such a neighborhood.

Consider first the case where $\operatorname{span}\{X_1, X_2\} \equiv \operatorname{span}\{\partial_x, x^{m-1} \partial_t\}$, in which microlocal analytic hypoellipticity is to be proved by the same method as in Section 7. Then $\tilde{\Theta}(x, t) \equiv x^{m-1}$. Let $L_\tau$ be the operator defined by conjugating the transpose of $L$ as in Section 7. The ordinary differential operators $A_t$ of that section should now be replaced by

$$\begin{aligned}
A_t = &\Big[(a_{1,1} - ia_{2,1})\partial_x + (a_{1,2} - ia_{2,2})i\tau(-1 - 2i(t'-t))\tilde{\Theta}\Big] \\
&\circ \Big[(a_{1,1} + ia_{2,1})\partial_x + (a_{1,2} + ia_{2,2})i\tau(-1 - 2i(t'-t))\tilde{\Theta}\Big] + O(\partial_x, \tau\tilde{\Theta}, 1)
\end{aligned}$$

where the coefficients in the term $O(\partial_x, \tau\tilde{\Theta}, 1)$ are $C^\omega$ and are bounded uniformly in $t, x, t', \tau$. Here the coefficients $a_{i,j}$ are still functions of $(x, t)$; we do not freeze coefficients as was done in Section 4. The principal part of $A_t$ is $-YY^*$, where

$$Y = (a_{1,1} - ia_{2,1})\partial_x - (a_{1,2} - ia_{2,2})i\tau\tilde{\Theta}$$

and for any operator $T$, we now denote by $T^*$ the formal adjoint of $T$, with respect to the usual Hilbert space structure of $L^2(\mathbb{R})$, rather than the transpose. The coefficients $a_{i,j}$ are real valued for $(x, t) \in \mathbb{R}^2$. Therefore for $(x, t) \in \mathbb{R} \times \mathbb{C}$,

$$-A_t = Y \circ Y^* + O(\partial_x, \tau\tilde{\Theta}, 1) + O(|\operatorname{Im}(t)|)O(\partial_x, \tau\tilde{\Theta})^2,$$

that is, $-A_t - YY^*$ is a quadratic polynomial, with real analytic coefficients, in $\partial_x, \tau\tilde{\Theta}, 1$ where $1$ denotes the identity operator, and the coefficients of all monomials of degree $2$ in $\partial_x, \tau\tilde{\Theta}$ are $O(|\operatorname{Im}(t)|)$.

For any $\phi \in C_0^2(\mathbb{R})$,

$$\int_{\mathbb{R}} YY^*\phi \cdot \bar{\phi} = \|Y^*\phi\|_{L^2(\mathbb{R})}^2 = \|Y\phi\|_{L^2(\mathbb{R})}^2 + \int_{\mathbb{R}} [Y, Y^*]\phi \cdot \bar{\phi},$$

and

$$[Y, Y^*] = -2\tau \operatorname{Det}(a)(\partial_x\tilde{\Theta}) + O(\partial_x, \tau\tilde{\Theta}, 1),$$

where norms without subscripts are $L^2(\mathbb{R}, dx)$ norms. By the pseudoconvexity hypothesis and the sign conventions $\operatorname{Det}(a) < 0$ and $\partial_x\tilde{\Theta} \geq 0$, the right hand side equals a nonnegative function, modulo an operator that is $O(\partial_x, \tau\tilde{\Theta}, 1)$. Therefore

$$2\int_{\mathbb{R}} YY^*\phi \cdot \bar{\phi} \geq \|Y\phi\|^2 + \|Y^*\phi\|^2 - C\int_{\mathbb{R}} |\phi|(|\partial_x\phi| + \tau|x^{m-1}\phi| + |\phi|),$$

which for large $\tau$ is

$$\geq c\int_{\mathbb{R}} |\partial_x\phi|^2 + \tau^{2/m}\left(1 + (\tau^{1/m}|x|)\right)^{2(m-1)}|\phi|^2 - C\int_{\mathbb{R}} |\phi|(|\partial_x\phi| + \tau|x^{m-1}\phi| + |\phi|),$$

by Cauchy-Schwarz. Replacing $YY^*$ by $-A_t$ introduces additional terms, but all these are majorized by

$$\varepsilon\|\phi\|_{\mathcal{H}_{0,\tau}^2}^2 + \varepsilon\int_{\mathbb{R}} |\partial_x\phi|^2 + \tau^{2/m}\left(1 + (\tau^{1/m}|x|)\right)^{2(m-1)}|\phi|^2$$
$$+ C\int_{\mathbb{R}} |\phi|(|\partial_x\phi| + \tau|x^{m-1}\phi| + |\phi|)$$

where $\varepsilon$ may be made as small as desired by taking $t, t'$ to be sufficiently small.

For $\tau > 0$ sufficiently large and $t, t'$ sufficiently small we thus obtain

$$-\operatorname{Re}\int_{\mathbb{R}} A_t\phi \cdot \bar{\phi} \geq c\int_{\mathbb{R}} |\partial_x\phi|^2 + \tau^2 x^{2(m-1)}|\phi|^2 - C\int_{\mathbb{R}} |\phi|(|\partial_x\phi| + \tau|x^{m-1}\phi| + |\phi|).$$

As in the proof of Lemma 7.2, this leads to the conclusion that for all sufficiently large $\tau$ and for all complex $t, t'$ and real $\rho$ sufficiently close to 0, $A_t$ maps $\mathcal{H}^2_{\rho,\tau}(\mathbb{R})$ invertibly to $\mathcal{H}^0_{\rho,\tau}(\mathbb{R})$, uniformly in all parameters. Reasoning just as in the proofs of Lemma 7.1 and Theorem 1.2 establishes the microlocal analytic hypoellipticity of $L$ in $\Gamma^+$, when $\tilde{\Theta}$ is divisible by $x^{m-1}$. $\qquad\square$

Consider next the case where $\tilde{\Theta}(0) = 0$ and $\tilde{\Theta}$ is not divisible by $x^{m-1}$, where it is to be shown that $L$ cannot be microlocally analytic hypoelliptic. The first obstacle to imitating the proof of Theorem 1.4 is that under our hypotheses, $(X_1 + iX_2)$ is not locally solvable near 0, whence neither is $(X_1 + iX_2)(X_1 - iX_2)$, so we cannot hope to solve the equation $LG_\lambda = F_\lambda$ for arbitrary $F_\lambda$ as was possible in Lemma 3.1. As a substitute, fix a bounded $C^\infty$ auxiliary function $n(\xi, \tau)$ which (i) is supported in the conic neighborhood $\{\tau < 0, \ |xi| < -\tau\}$ of the ray $\{\xi = 0, \tau < 0\}$, (ii) is identically equal to 1 in a smaller conic neighborhood minus some compact set, and whose inverse Fourier transform (iii) is $C^\omega$ except at the origin, and (iv) agrees with a Schwartz function away from the origin.[7] Define the Fourier multiplier operator $(Bf)\hat{\ }(\xi, \tau) = n(\xi, \tau)\hat{f}(\xi, \tau)$.

**Lemma 8.2.** *Let $\{X_1, X_2\}$ satisfy the bracket and pseudoconvexity hypotheses. Then there exists a neighborhood $U$ of 0 such that for all $f \in L^2(U)$ there exist $v, g \in L^2(U)$ satisfying $Lg + Bv = f$ in $U$, such that the $L^2(U)$ norms of $g, v$ are bounded by a fixed constant times the norm of $f$.*

*Proof (sketch).* $L$ is elliptic outside $\Sigma$ and is subelliptic in a conic neighborhood of $\Sigma^+$, because of the bracket and pseudoconvexity hypotheses, as follows from the reasoning of Kohn [14]. On the other hand, the pseudodifferential operator $B$ is elliptic in a conic neighborhood of $\Sigma^-$. This is all that is needed for the proof of the nearly identical Lemma 7.2 of [9] to apply. $\qquad\square$

Defining $F_\lambda$ as in Lemma 3.1, solve $LG_\lambda = F_\lambda + Bv_\lambda$ with $v_\lambda = O(1)$ in $L^2$ norm. Proceeding as in Sections 3 and 4 with the notations and definitions introduced there, $e^{-i\lambda^m t}G_\lambda$ satisfies a conjugated equation, with right hand side $e^{-i\lambda^m t}F_\lambda + e^{-i\lambda^m t}Bv_\lambda$. We are interested in this equation in the region where $\mathrm{Im}(t) < 0$ and $|\mathrm{Im}(t)| = \lambda^{-1/q}|\mathrm{Im}(s)| \sim \lambda^{-1/q}$, and hence $|\exp(-i\lambda^m t)| = O(\exp(-c\lambda^{m-q^{-1}}))$ for some $c > 0$. Because $Bv_\lambda \in L^2(U)$ uniformly in $\lambda$ and because the Fourier transform of $Bv_\lambda$ is supported where $\tau < 0$ and $|\xi| < |\tau|$, $Bv_\lambda$ extends to a holomorphic function of $t$ near 0 in the half space $\{\mathrm{Im}(t) < 0\}$, and each partial derivative of $Bv_\lambda$, with respect to $(x, t)$, is $O(\lambda^N)$ there for some finite $N$; this follows by writing $Bv_\lambda$ as the inverse Fourier transform of its Fourier transform, and noting that the resulting integral converges absolutely for $\text{Im}(t) < 0$. Therefore

$$e^{-i\lambda^m t} Bv_\lambda = O(\exp(-c\lambda^{m-q^{-1}})) = O(\exp(-c\lambda^p))$$

in this region, for some $c > 0$, as $\lambda \to +\infty$.

---

[7]Such functions $n$ may be constructed using the cutoff functions of Andersson and Hörmander, which are discussed in [13].

The situation is therefore the same as in Sections 3 and 4, except for the presence of an extra term of magnitude $O(\exp(-c\lambda^p))$ on the right hand side of the conjugated equation for $\exp(-i\lambda^m t)G_\lambda(x, t)$. Such a term has no effect in the argument immediately following Lemma 4.2, Following the nor in the proof of Lemma 4.2 outlined in Section 9. In other respects the analysis of $L = (X_1 + iX_2)(X_1 - iX_2) + c_1 X_1 + c_2 X_2 + c_3$ is identical to that for $X_1^2 + X_2^2$, modulo minor changes in formulas. In this case the existence of eigenvalues was established for all $Q$ in Theorem 5.2, so it follows in full generality that $L$ cannot be analytic hypoelliptic if $\tilde{\Theta}$ is not divisible by $x^{m-1}$.

## 9. Proof of Lemma 4.2

In this section we outline the proof of the remaining step in the reduction of analytic hypoellipticity, for sums of squares of two vector fields in $\mathbb{R}^2$, to eigenvalue problems for ordinary differential operators. The proof is a straightforward adaptation of that given for the corresponding Lemma 4.3 of [3], with certain systematic changes in formulas. We will therefore give the definitions, notations, and statements of sublemmas, armed with which the determined reader will be able to construct a full proof by following the exposition of [3] line by line.

Recall the conjugated operator

$$B_\lambda = P\Big(\lambda^{-1} y, \lambda^{-1/q}(\zeta_0 + re^{i\theta}), \partial_y, i[\Theta + R_\lambda](1 - i\lambda^{-p}\partial_s)\Big),$$

where $p = m - q^{-1}$ and $q < \infty$. Assume that

$$B_\lambda u_\lambda = \psi + O(e^{-\delta\lambda^p})$$

for all large $\lambda \in \mathbb{R}^+$. The extra term $O(\exp(-\delta\lambda^p))$ on the right is absent in Lemma 4.2, but occurs in the analogue needed for the proof of Theorem 1.8. This equation holds for all $(y, s) \in \mathbb{R} \times \mathbb{C}$ satifying $|y| \leq c\lambda$ and $C^{-1} \leq -\text{Im}(s) \leq C$, for some $c > 0$ and for any $C \in \mathbb{R}^+$, where $\delta$ depends on $C$. By (4.1) we know that $u_\lambda$ and all its partial derivatives of order $\leq 2$ are bounded by $\exp(B\lambda^p)$, for some large constant $B$, for all $(y, s)$ in this same region. What is desired is essentially a bound $u_\lambda = O(\exp(\mu\lambda^p))$, for arbitrarily small $\mu > 0$, in the smaller region where $|y| \leq \lambda^{1-\varepsilon_0}$ and $s$ belongs to a neighborhood of $\zeta_0$.

The proof requires the selection of auxiliary constants $A_0, \nu, N, \sigma_0, \ldots \sigma_N, A, \gamma$ satisfying various constraints. These depend on $B, \mu, \zeta_0$, but not on $\lambda$. Let $A_0$ be a large constant to be chosen later, and set $C_0 = 2|\zeta_0|$.

**Lemma 9.1.** *There exist* $0 < \nu < \mu$, $N < \infty$, *and* $3C_0 > \sigma_0 > \sigma_1 > \cdots > \sigma_N > 2C_0$ *such that*

$$\tfrac{1}{2}\mu < A_0 - (N+1)\nu < \mu$$

*and*

$$(s - i\sigma_j)(1+r)^{1/qm} \notin E^{1/q} \tag{9.1}$$

*for all* $s \in \mathbb{R}$, $0 \leq j \leq N$ *and* $r \in \mathbb{R}^+$ *satisfying*

$$(A_0 - j\nu) \geq r \geq (A_0 - (j+1)\nu).$$

A small twist needed to adapt the proof of the corresponding lemma in [3] may be found in [21].

The constants $\gamma, A$ will be required to satisfy finitely many constraints, to be encountered in the course of the proof, all of the two forms

$$\gamma \text{ and } \gamma \cdot A \text{ are sufficiently small} \tag{9.2}$$

$$\gamma \cdot A^2 \text{ is sufficiently large} \tag{9.3}$$

relative to various quantities depending on $B, \mu, \zeta_0, \nu, N, \sigma_0, \ldots \sigma_N$. Any finite collection of such constraints is satisfied by some pair $\gamma, A \in \mathbb{R}^+$; we assume henceforth that $\gamma, A$ do satisfy them.

The first constraint on $\gamma, A$ is that

$$(s - i\sigma_j)\big[1 + r - i2\gamma(s - i\sigma_j)\big]^{1/mq} \notin E^{1/q} \tag{9.4}$$

for all $s \in \mathbb{R}$ satisfying $|s| \leq 4A$, all $0 \leq j \leq N$, and all $r \in I_j$ where

$$I_j = [A_0 - (j+1)\nu, A_0 - j\nu].$$

This follows from (9.1), using Lemma 3.7, as shown in [3].

Extend the coefficients of $B_\lambda$ so as to be globally defined with respect to $y$, as was done in the proof of Lemma 4.1. Define

$$f_\sigma(y, s) = e^{-\gamma\lambda^p(s - i\sigma)^2} u_\lambda(y, s - i\sigma)\eta(A^{-1}s)$$

where $\eta \in C_0^\infty(\mathbb{R})$ is supported in $(-4, 4)$ and is identically equal to 1 on $[-2, 2]$. Set

$$\mathbb{L}_\sigma = e^{-\gamma\lambda^p(s - i\sigma)^2} \circ B_\lambda \circ e^{\gamma\lambda^p(s - i\sigma)^2}$$

$$= P\Big(\lambda^{-1}y, \lambda^{-1/q}(\zeta_0 + re^{i\theta}), \partial_y, i[\Theta + R_\lambda](1 - 2i\gamma(s - i\sigma) - i\lambda^{-p}\partial_s)\Big),$$

Then

$$\mathbb{L}_\sigma f_\sigma = \Psi_\sigma \quad \text{where } \Psi_\sigma(y,s) = e^{-\gamma\lambda^p(s-i\sigma)^2}\psi(y) + O(e^{-\delta\lambda^p})e^{-\gamma\lambda^p(s-i\sigma)^2}$$

for all $s \in \mathbb{R}$ satisfying $|s| \leq 2A$, provided that $\sigma \geq C^{-1}$. For $|s| > 2A$ the cutoff function $\gamma$ comes into play, but the factor $\exp(-\gamma\lambda^p(s-i\sigma)^2)$ is $O(\exp(\lambda^p[C-\gamma A^2]))$, with $C$ dependent only on quantities fixed before $\gamma, A$ are chosen. Thus both $\Psi_\sigma$ and $f_\sigma$ are $O(\exp(-C_2\lambda^p))$ where $|s| \geq A/2$, for any desired $C_2 < \infty$, provided that $\gamma A^2$ is chosen to be sufficiently large. The same holds for their partial derivatives of first and second orders with respect to $y, s$.

Define the partial Fourier transform

$$\hat{f}_\sigma(y,\xi) = \int_{\mathbb{R}} f_\sigma(y,s)e^{-is\xi}\,ds$$

for all $\xi \in \mathbb{R}$.

**Lemma 9.2.** *For any constant $A_0 < \infty$, for all sufficiently large $\lambda$, for all $|\xi| \leq A_0\lambda^p$, $|y| \leq \lambda^{1-\varepsilon_0}$ and $0 \leq \sigma, \sigma' \leq 3C_0$,*

$$\hat{f}_\sigma(y,\xi) = e^{(\sigma-\sigma')\xi}\hat{f}_{\sigma'}(y,\xi) + O(e^{-\lambda^p}).$$

This is proved by shifting the contour of integration with respect to $s$ in the definition of the partial Fourier transform. For $|s| \geq A$ the integrand is not holomorphic, resulting in an error term that is $O(\exp(-\delta\gamma A^2\lambda^p + C\lambda^p))$ for some $\delta > 0$, which is negligible if $\gamma A^2$ is sufficiently large.

**Lemma 9.3.** *There exists a large constant $A_0$ such that for all sufficiently large $\lambda$ and all $|y| \leq \lambda^{1-\varepsilon_0}$, for all $|\operatorname{Im}(\zeta_0)|/2 \leq \sigma \leq 3C_0$,*

$$\int_{|\xi| \geq A_0\lambda^p} \left(|\hat{f}_\sigma(y,\xi)|^2 + |\partial_y\hat{f}_\sigma(y,\xi)|^2\right)d\xi \leq e^{-2\lambda^p}.$$

This is proved by combining the preceding lemma with the weak bound $O(\exp(B\lambda^p))$ for $u_\lambda$, choosing $A_0$ sufficiently large to overwhelm the factor $\exp(B\lambda^p)$.

For functions $\varphi$ of $y \in [-\lambda^{1-\varepsilon_0}, \lambda^{1-\varepsilon_0}]$ define the dual norms

$$\|\varphi\|_{H_\rho}^2 = \int_{|y| \leq \lambda^{1-\varepsilon_0}} |\varphi(y)|^2 \langle y \rangle^{-2(m-1)} e^{+\rho|y|^m}\,dy,$$

$$\|\varphi\|_{H_\rho^*}^2 = \int_{|y| \leq \lambda^{1-\varepsilon_0}} |\varphi(y)|^2 \langle y \rangle^{+2(m-1)} e^{-\rho|y|^m}\,dy$$

where $\rho > 0$ is fixed and small. Denote by $\hat{f}(\xi)$ the function $y \mapsto \hat{f}(y,\xi)$.

**Lemma 9.4.** *There exists $\delta > 0$ such that for all sufficiently large $\lambda$, for all $0 \leq j \leq N$, whenever $\lambda^{-p}\xi \in I_j$,*

$$\|\hat{f}_{\sigma_j}(\xi)\|_{H^*_\rho}^2 \leq e^{-\delta\lambda^p} \int_{|\eta| \leq A_0\lambda^p} \|\hat{f}_{\sigma_j}(\eta)\|_{H^*_\rho}^2 \, d\eta + e^{-\delta\lambda^p}.$$

Lemma 4.2 follows from this result and Lemmas 9.3 and 9.2 by an induction on $j$ as in [3], so it remains only to prove Lemma 9.4. This requires a final sublemma.

Denote by $\mathbb{L}^*_\sigma$ the formal transpose of $\mathbb{L}_\sigma$ with respect to the pairing

$$\langle f, \, g \rangle = \iint_{\substack{|y| \leq \lambda^{1-\varepsilon_0} \\ |s| \leq A}} f(y,s)g(y,s) \, dy \, ds,$$

ignoring all boundary terms arising from integration by parts in the formal identity $\langle \mathbb{L}_\sigma f, \, g \rangle = \langle f, \, \mathbb{L}^*_\sigma g \rangle$.

**Lemma 9.5.** *There exist $\delta > 0$ and an open set $\mathbb{C} \supset \Omega \supset \{s \in \mathbb{R} : |s| \leq A\}$ such that for all sufficiently large $\lambda$ and sufficiently small $\rho > 0$, for all $0 \leq j \leq N$, $\xi \in \lambda^{-p}I_j$ and $\varphi \in L^2(\mathbb{R})$ supported in $[-\lambda^{-\varepsilon_0}, \lambda^{-\varepsilon_0}]$, there exist $g, E$ defined on $[-\lambda^{1-\varepsilon_0}, \lambda^{1-\varepsilon_0}] \times \Omega$ satisfying*

$$e^{is\xi}\mathbb{L}^*_{\sigma_j}(e^{-is\xi}g)(y,s) = \varphi(y) + E(y,s) \qquad in \; [-\lambda^{1-\varepsilon_0}, \lambda^{1-\varepsilon_0}] \times [-A,A],$$

$$\|g(\cdot,s)\|_{\mathcal{H}^2_\rho} \leq C\|\varphi\|_{\mathcal{H}^0_\rho} \qquad for \; all \; s \in \Omega$$

$$\|E(\cdot,s)\|_{\mathcal{H}^0_\rho} \leq Ce^{-\delta\lambda^p}\|\varphi\|_{\mathcal{H}^0_\rho} \qquad for \; all \; s \in \Omega.$$

The norms $\mathcal{H}^k_\rho$ are those of $\mathcal{H}^k_\rho([-\lambda^{1-\varepsilon_0}, \lambda^{1-\varepsilon_0}])$, defined as in Section 3 except that the integration with respect to $y$ extends only over $[-\lambda^{1-\varepsilon_0}, \lambda^{1-\varepsilon_0}]$, rather than over the entire real line. The proof of this lemma is very similar to that of Lemma 4.1; the relation (9.4) together with Lemma 3.6 guarantee that those ordinary differential operators arising in the analysis are indeed invertible. Lemma 9.4 is then deduced by a duality argument like that used in Section 4 to deduce Theorem 1.4 from Lemma 4.1, shifting the contour of integration with respect to $s$ to majorize the principal term resulting after integration by parts, as in the proof of Lemma 6.4 of [3].

## 10. COMMENTS

1. Denote by $\Sigma$ the characteristic variety of $L$, that is, the set of all points in the cotangent bundle where the principal symbol of $L$ vanishes. In our analysis there is no essential distinction between the four coefficients $\tilde{\Theta}$:

- $x(x^2 + t^2)$, for which $\Sigma$ is a (smooth) symplectic manifold,
- $x(x^2 - t^2)$, for which $\Sigma$ is a singular variety with normal crossings,

- $x(x^2 - t^3)$, for which $\Sigma$ is more singular, and
- $x^2 + t^2$, for which $\Sigma$ consists of the single fiber $\{x = t = \xi = 0\}$,

and no distinction whatsoever between the first two examples. Yet in terms of symplectic geometry, the four characteristic varieties have little in common.[8]

2. One of the principal hypotheses in the work of Tartakoff [18],[19] and of Treves [20] giving sufficient conditions for analytic hypoellipticity is that the characteristic variety be symplectic. Treves had conjectured that analytic hypoellipticity should fail to hold (for, say, any $X_1^2 + X_2^2$) whenever the characteristic variety $\Sigma$ contains a smooth nonconstant curve $\gamma$ that is orthogonal to the tangent space of $\Sigma$ (at each point of $\gamma$). For the case of two linearly independent vector fields in $\mathbb{R}^3$, any null bicharacteristic for some nonvanishing vector field in the span of $X_1, X_2$ that lies above the set of points in $\mathbb{R}^3$ where $X_1, X_2, [X_1, X_2]$ are linearly dependent is such a curve $\gamma$, and the conjecture has been established [7] in this case. More recently, the conjecture has also been confirmed for a three dimensional example in which the only such curve $\gamma$ is a fiber of the cotangent bundle [21].

Himonas and Treves asked for examples of operators not analytic hypoelliptic, yet not possessing any such curves $\gamma$. The history of this question is a bit obscure; Oleĭnik and Radkevič had already published by 1973 [16] examples of analytic nonhypoelliptic operators with characteristic varieties $\Sigma$ that are symplectic and hence possess no such curves $\gamma$, but had not explicitly discussed the symplectic geometric nature of $\Sigma$. Much later it was observed that our examples $\partial_x^2 + (x^k(x^2 + y^2)\partial_t)^2$ likewise have symplectic characteristic varieties $\{x = \xi = 0\}$, and by Proposition 1.7 are not analytic hypoelliptic for $k = 1$ or $k$ even. Special cases of the examples of Oleĭnik and Radkevič were subsequently rediscovered by Hanges and Himonas.[9]

3. Our results should generalize to sums of squares of any finite number of vector fields in $\mathbb{R}^2$. The necessary and sufficient condition for analytic hypoellipticity should be the existence of a coordinate system in which span$\{X_j\} \equiv$ span$\{\partial_x, x^{m-1}\partial_t\}$, that is, both $\partial_x$ and $x^{m-1}\partial_t$ may be expressed as linear combinations, with $C^\omega$ coefficients, of the $X_j$, and conversely.

4. In Proposition 1.7 we have singled out certain examples for which Conjecture 1.5 can be proved relatively easily. With more labor, similar reasoning should apply to other special cases, as well, but the conjecture can be formulated more generally by allowing the coefficients of $Q$ to be arbitrary complex numbers, and then neither the *ad hoc* proof of Proposition 1.7 nor the method of proof of Theorem 5.2 appear likely to lead to a proof in full generality.

---

[8]This accords with a remark of Métivier [15]: "Il ne faut sans doute pas vouloir relier trop rigidement le condition "$\Sigma$ symplectique" à l'hypoellipticité analytique."

[9]A higher order example is $L = L_1 \circ L_2$, where the $L_j$ are chosen to have characteristic varieties $\Sigma_j$ such that $\Sigma_2 \subset \Sigma_1$, $\Sigma_1$ is symplectic, and $L_2$ is not analytic hypoelliptic. Then $L$ is not analytic hypoelliptic, and its characteristic variety is $\Sigma_1$. An explicit example is $(\partial_x^2 + x^{2k}\partial_t^2) \circ (\partial_x^2 + (x^2 + t^2)\partial_t^2)$.

A more promising approach is related to the phenomenon of Stokes lines. A solution $\psi_z^+$ having the required asymptotics as $x \to +\infty$ has this behavior in a conic neighborhood of $\mathbb{R}^+$, but not for $x$ in all of $\mathbb{C}$; the Stokes lines separate disjoint sectors in the complex plane in which the solution has distinct asymptotics. Whenever a (normalized, homogeneous) polynomial $Q(x, z)$ is not identically equal to $x^{m-1}$, then $x \mapsto Q(x, z)$ has for each $z$ at least two distinct zeros in the complex plane. As in the proof of Theorem 5.2, these zeros should govern the asymptotic behavior of $W(z)$ as $|z| \to \infty$. One could attempt to show that there exist zeros $w_i(z)$ and disjoint sectors $\Gamma_i$, for $i = 1, 2$, such that the asymptotic behavior of the Wronskian $W(z)$ as $z$ tends to infinity through $\Gamma_i$ is governed by $w_i(z)$; more precisely, that for $z \in \Gamma_i$, $W(z) \sim \exp(c_i z^m) z^{\gamma_i} (1 + O(|z|^{-\delta})$ for some exponents $\gamma_i$. If one could show further that $c_1 \neq c_2$, then $W$ could not have a polynomial logarithm and hence would have zeros. One advantage of this approach is that it requires no analysis of the exponents $\gamma_i$ (let alone of higher order terms in the asymptotic expansions). Some of the ideas of Yu [23] should be useful here.

5. The fact that ordinary differential operators play such a prominent role in our analysis is a feature of the particular class of low-dimensional partial differential operators under consideration. Globally elliptic operators in more than one variable play the corresponding role in other situations.

6. In its present form, the method employed here is too primitive to apply to the general case of two independent vector fields in $\mathbb{R}^3$.

MICHAEL CHRIST, DEPARTMENT OF MATHEMATICS, UNIVERSITY OF CALIFORNIA, LOS ANGELES, CA 90095-1555

*E-mail address*: `christ@math.ucla.edu`

*Current address*: Mathematical Sciences Research Institute, 1000 Centennial Road, Berkeley, CA 94720